\documentclass[journal]{IEEEtran}
\IEEEoverridecommandlockouts                              
\usepackage{setspace}
\usepackage{graphics} 
\usepackage{epsfig} 
\usepackage{mathptmx} 
\usepackage{times} 
\usepackage{amsmath} 
\usepackage{amssymb}  
\usepackage{epstopdf}
\usepackage{cite}
\usepackage[final]{changes}
\usepackage{xcolor}
\usepackage{bbm}
\usepackage{tikz}
\usetikzlibrary{arrows}
\usepackage{caption}
\usepackage{subcaption}
\usepackage{algpseudocode}
\usepackage{mathtools}
\usepackage{algorithm}
\usepackage{balance}
\usepackage{multirow}
\usepackage{bm}
\allowdisplaybreaks

\newtheorem{theorem}{Theorem}[section]

\newtheorem{proposition}[theorem]{Proposition}
\newtheorem{problem}[theorem]{Problem}
\newtheorem{corollary}[theorem]{Corollary}
\newtheorem{definition}[theorem]{Definition}
\newtheorem{assumption}[theorem]{Assumption}
\newtheorem{remark}[theorem]{Remark}

\ifCLASSINFOpdf

\else

\fi

\begin{document}
%
\title{\LARGE \bf Optimal Co-design of Industrial Networked Control Systems with State-dependent Correlated Fading Channels}

\author{Bin Hu 
	\thanks{Bin Hu is with Department of Engineering Technology, Old Dominion University, Norfolk, VA 23529, USA.
		{\tt\small bhu@odu.edu}}  
	and Tua A. Tamba
	\thanks{T. A. Tamba is with Department of Electrical Engineering (Mechatronics), Parahyangan Catholic University, Bandung 40141, West Java, Indonesia
		{\tt\small ttamba@alumni.nd.edu}}%
}

\maketitle

\begin{abstract}
This paper examines a co-design problem for industrial networked control systems (NCS) whereby physical systems are controlled over wireless fading channels.
In particular, the considered wireless channels are also stochastically dependent on the physical states of moving machineries  in an industrial working space.
In this paper, the moving machineries are modeled as Markov decision processes whereas the characteristics of the correlated fading channels are modeled as a binary random process whose probability measure is a function of both the moving machineries' physical states and the channels' transmission power.
Under such state-dependent fading channel models, sufficient conditions which ensure the stochastic safety of the NCS are first derived.
Using the derived safety conditions, the co-design problem is then formulated as a constrained joint optimization problem that seeks for  optimal control and transmission power policies which simultaneously minimize both the communication and control costs in an infinite time horizon.
This paper shows that such optimal co-design policies can be obtained in an efficient manner from the solution of convex programs.
Simulation results from industrial NCS models consisting forklift truck and networked DC motor systems are also presented to illustrate and verify both the advantages and efficacy of the proposed co-design solution framework.
\end{abstract}

	\renewcommand{\abstractname}{Note to Practitioners}
\begin{abstract}
This paper is motivated by the problems of designing efficient communication and control policies to ensure the safety of factory automation where different processes coordinate with each other through wireless networks. One of the main challenges for this problem lies in the fact that wireless networks used for safety are highly unreliable, and can be seriously disrupted by operational machinery in the vicinity. Existing approaches that decouple the design of communication and control policies may fail to achieve efficiency for factory automation systems due to the interaction between the communication~(cyber) and physical systems. By taking into account such \emph{cyber-physical} couplings, this paper develops a novel co-design framework under which the communication and control policies are coordinated to achieve both system safety and efficiency.  Under the co-design framework, this paper further shows that the communication and control policies that minimize the use of both communication and control resources in the long run while respecting safe operations, can be computed efficiently. This allows the proposed co-design method to go beyond the simple example illustrated in this paper and apply to more complex practical systems, such as automobile assembly system, manufacturing factory with automated heavy facilities, and automated warehouse with mobile industrial robots. 
\end{abstract}

\IEEEpeerreviewmaketitle

\section{Introduction}
\label{sec: introduction}
\subsection{Background and Motivation}
Over the last few decades, wireless communication technologies have rapidly evolved and continuously developed to support and improve various industrial processes automation.
Modern industrial automation architectures are now often equipped with industrial wireless communication protocols such as WirelessHart \cite{song2008wirelesshart} and WiMAX \cite{gungor2009industrial} which offer the promise of process  safety and efficiency  improvements.
Compared to more traditional industrial automation systems which mainly rely on expensive wired communication systems,  wireless communication technologies have been known to be cost-effective and enables more flexible, intelligent and productive operations of industrial automation systems.
However, conventional industrial wireless communication protocols are also known to be inherently unreliable and often  subject to channel fading phenomena \cite{agrawal2014long,aakerberg2011future,willig2005wireless,willig2008recent}. 
For instance, recent works in \cite{agrawal2014long,kashiwagi2010time,quevedo2013state,olofsson2016modeling} have shown that the characteristics of channel fading in industrial environments are often statistically dependent on the physical motions of large machineries which operate in industrial environment. 
In particular, such state-dependent fading channels correlate physical states of the industrial systems with the states of the communication channel, thereby introducing great challenges on the need for assuring both the safety and efficiency of the entire industrial operations.

This paper examines an optimal co-design problem for industrial NCS in the presence of state dependent correlated fading channels. 
Specifically, this paper examines the so-called shadow fading phenomenon which may occurs from temporary obstruction of the used radio signals transmissions by the movements of large and heavy machineries. 
Such a fading thus represents a type of communication channel failure which may compromises the overall safety of the industrial NCS.

From the view point of wireless communication systems design, one may in principle alleviates the channel fading effect by increasing the transmission power of the communication system \cite{tse2005fundamentals}. 
However, such an increase may in the long run leads to overuse of energies and consequently compromises the overall system efficiency.
Motivated by this trade off between system safety and efficiency,  this paper examines a co-design problem in industrial NCS. 
This problem essentially concerns with the development of control policy for the moving machineries on the one hand and communication policy for network system protocol on the other hand but simultaneously able to achieve optimal control and communication performances while at the same time satisfy the safety specifications. 
It is now well-understood that the main technical challenge in solving such a problem lies in the strong coupling between control and communication policies that is induced by the presence of state-dependent fading phenomena on the communication channels. 
This paper contributes to the development of solution methods for this problem by adopting a novel state-dependent channel model from \cite{hu2017co} to formulate  linear/quadratic program for finding the optimal policies.


\subsection{Relevant Work and Contributions}

The co-design issue addressed in this paper is one of the challenges in NCS development with regard to the limited capacity or imperfection of communication channels.
Out of concerns on performance degradation that may be caused bysuch an imperfection, this issue has recently attracted a great amount of interests from researchers in control and communication communities.
It is beyond the scope of this paper to exhaustively review all of the present literature on this issue. 
Instead, this section of the paper will be focused on reviewing recent works on co-design problems in NCS with an emphasis on their industrial applications.

The formulation of co-design problems in NCS are often reduced to joint stabilization or optimization problems under the constraints of limited channel capacity and control resources.
These joint stabilization/optimization problems may be classified into several categories according to the characteristics of the resources/channels (e.g. transmission power \cite{gatsis2014optimal, ren2018infinite, wen2018transmission,hu2017co}, data rate \cite{rabi2016separated, leong2017event}, bandwidth \cite{peng2013event, peng2013novel} or network topology \cite{zhao2008integrated}) as well as the imperfections (e.g. delay and packet dropout \cite{gatsis2014optimal, hu2017co,leong2017event}) of the used communication systems.
With the goal of ensuring the energy efficiency of the NCS, the work in \cite{gatsis2014optimal} has shown that the optimal co-design of both transmission power and control policies for a linear dynamical system over a wireless fading channel may somewhat be  decoupled if a restricted set of information structure is assumed. 
Similar idea and results were reported in \cite{ren2018infinite} where optimal power and remote estimation policies were jointly designed to minimize an infinite horizon cost which consists of power utilization and estimation errors. 
Furthermore, recent results reported in \cite{varma2016energy} has  show that the separation principle introduced in \cite{ren2018infinite, gatsis2014optimal} may also be applied to more general system structures whenever certain time-triggering conditions are satisfied.

From the viewpoint of communication bandwidth utilization, the works in \cite{peng2013event, peng2013novel,leong2017event,varma2016energy,lyu2017co,zhao2008integrated,peters2016controller} have examined related stability and optimality issuesin the co-design problem of NCS.  
By adopting an event-triggered communication scheme, the works in \cite{peng2013event, peng2013novel} proposed a solution approach to simultaneously design an event-triggered communication scheme and feedback controllers that may assure the $\mathcal{H}_{\infty}$ \cite{peng2013event} or $\mathcal{L}_{2}$ performance \cite{peng2013novel} of the NCS. 
In \cite{leong2017event,zhao2008integrated}, the co-design problem was formulated as a joint optimization problem which searches for optimal transmission scheduling and control policies that minimize both the control and communication costs. 
Under an assumption that the wireless fading channels are independent of the physical system, the work in \cite{leong2017event} showed that the joint optimization problem can be reduced to a separation design where the optimal controller is defined by the solution of an LQG type problem whereas the optimal scheduling policy is triggered by the estimation error covariance at the controller. 


All of the aforementioned results essentially assumed that the used communication channels  are stochastically independent of the physical systems. 
Such an assumption, however, generally does not hold in most industrial environments \cite{kashiwagi2010time, quevedo2013state,olofsson2016modeling,agrawal2014long,eriksson2016long} in which channel fading are often highly correlated with the states of the operating machineries.
Such a high correlation property consequently invalidates the use of the separation principle to solve the problem of joint design of optimal control and communication policies. 
Motivated by this realization, the work in \cite{quevedo2013state} then explored the state estimation problem of  linear NCS with correlated fading channels. 
The results in \cite{quevedo2013state} were further extended to address the co-design of optimal power control and coding scheme \cite{quevedo2014power} as well as network reconfiguration under correlated fading channel \cite{leong2016network}.


The co-design framework discussed in this paper is different from those in \cite{quevedo2013state, quevedo2014power}  since the former examines stochastic safety and efficiency of industrial NCS under state-dependent correlated fading channels. 
Specifically, the NCS in the proposed framework is controlled over a wireless fading channel whose states are correlated with the states of the existing moving industrial machineries (modeled by a Markov decision process). 
As such, this paper formulates the co-design problem as a constrained joint optimization problem whereby the guarantee for stochastic safety on the one side is enforced through some predefined safety constraints.  The assurance of stochastic efficiency on the other side is determined by the existence of optimal transmission power and control policies which minimize the expectation of an infinite horizon, joint communication and control cost.  Below, we summarize the main contributions of the present paper.
\begin{itemize}
	\item  This paper's first contribution is the derivation of sufficient conditions that assure the stochastic safety of the considered industrial NCS. Under newly defined mild assumptions, the derived conditions are different from those in \cite{hu2017co} since the former are convex (i.e. have linear or quadratic structures) whereas the latter were non-convex.  This mens that the newly proposed conditions can be examined computationally in a more efficient manner.  
	\item  The second contribution of this paper is the formulation of efficient computational methods for solving the co-design problem of NCS. Specifically, using newly derived safety sufficient conditions, this paper shows that the co-design problem can be solved by either linear programming (for the weak notion of safety in expectation) or quadratic programming (for the strong notion of almost sure asymptotic safety) method. 
	The formulations of such computational methods essentially define a significant difference between this paper and that of \cite{hu2017co}. 
\end{itemize}

The remainder of this paper is organized as follows. 
The description and modeling framework of the considered industrial NCS is presented in Section \ref{sec:system-framework}. 
Section \ref{sec: problem-formulation} formulates the co-design problem within the notions of stochastic safety and efficiency. 
Under some mild assumptions regarding the characteristics of the considered industrial NCS as defined in Section \ref{sec:assumptions},  sufficient conditions for guaranteeing the stochastic safety of the NCS are derived in Section \ref{sec: safety}. 
Using these safety conditions, Section \ref{subsec: co-design} presents a constrained two-player cooperative game framework and two convex programs (linear and quadratic programs) for solving the co-design problem of NCS. 
Simulation results that illustrate the applications of the proposed framework are given in Section \ref{sec: simulation}. Section \ref{sec: conclusion} concludes the paper with discussion and remarks.

\vspace{8pt}
\noindent
{\bf Notations:} Throughout the paper, the sets of non-negative real numbers and integers are denoted as $\mathbb{R}_{\geq 0}$ and $\mathbb{Z}_{\geq 0}$, respectively. 
An $n$ dimensional vector space is denoted by $\mathbb{R}^{n}$. 
The infinity norms of a vector $x \in \mathbb{R}^{N}$ and a matrix $A$ are denoted by $|x|$ and $\|A\|$, respectively. 
A function $f(k)$ is said to be essentially ultimately bounded if $\exists M > 0$ such that $ |f(k)|_{\mathcal{L}_{\infty}}=\text{ess}\sup_{k \geq 0} \|f(k)\|_2 \leq M$ where $\|\cdot\|_2$ is the Euclidean norm. 
A function $\alpha(\cdot): \mathbb{R}_{\geq 0} \rightarrow \mathbb{R}_{\geq 0}$ is said to be of class $\mathcal{K}$ if it is continuous, strictly increasing, and satisfies $\alpha(0)=0$. 
A function $\alpha(\cdot)$ is  said to be of class $\mathcal{K}_{\infty}$ if it is a class $\mathcal{K}$ function and radially unbounded. 
A function $\beta(\cdot, \cdot): \mathbb{R}_{\geq 0} \times \mathbb{R}_{\geq 0} \rightarrow \mathbb{R}_{\geq 0}$ is  said to be of class $\mathcal{KL}$ function if $\beta(\cdot, t)$ is a class $\mathcal{K}_{\infty}$ function for each fixed $t \in \mathbb{R}_{\geq 0}$ and in addition satisfies $\beta(s, t) \rightarrow 0$ for each $s \in \mathbb{R}_{\geq 0}$ as $t \rightarrow +\infty$. 
A function $\overline{\beta}(\cdot, \cdot, \cdot): \mathbb{R}_{\geq 0} \times \mathbb{R}_{\geq 0} \times \mathbb{R}_{\geq 0} \rightarrow \mathbb{R}_{\geq 0}$ is said to be of class $\mathcal{KLL}$~($\overline{\beta} \in \mathcal{KLL}$), if for each $r \geq 0$, $\overline{\beta}(\cdot, \cdot, r) \in \mathcal{KL}$ and $\overline{\beta}(\cdot, r, \cdot) \in \mathcal{KL}$ hold. 
The expected value and probability measure of a random variable $\{y\}$ are denoted as $\mathbb{E}\{y\}$ and $\mathbb{P}\{y\}$, respectively.

\section{Heterogeneous System Framework}
\label{sec:system-framework}

\begin{figure}[!b]
\centering
\includegraphics[scale=.25]{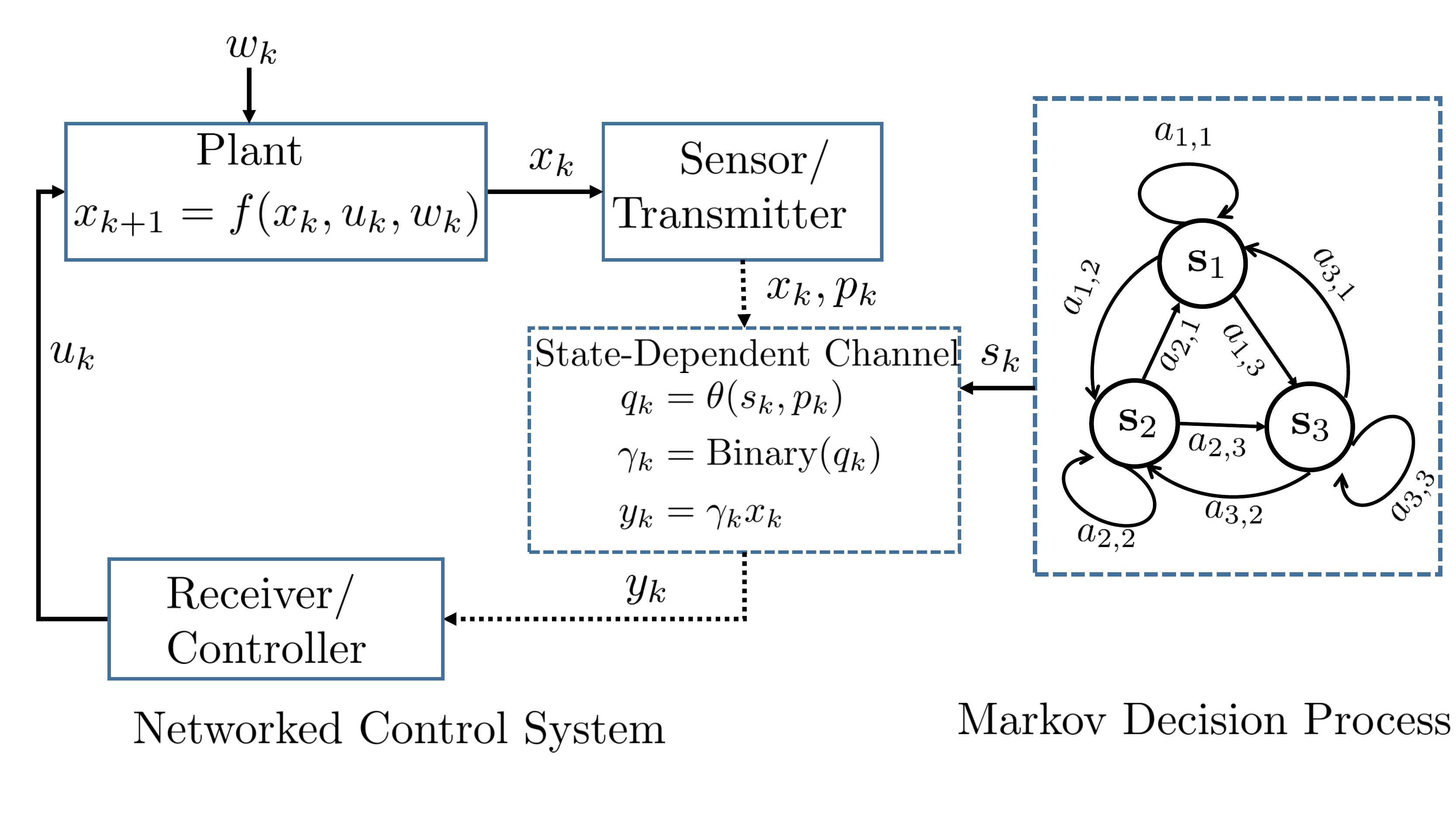}
\caption{Heterogeneous system framework with NCS and MDP.}
\label{fig: system-framework}
\end{figure}

Fig. \ref{fig: system-framework} illustrates a heterogeneous system modeling framework that is used in this paper. 
This figure shows that the framework essentially consists of two main subsystems. 
The first subsystem is a Markov decision process~(MDP) model describing a discrete event decision making process which manages high level control tasks. 
The second, lower level subsystem is a discrete time NCS model which consists of a nonlinear physical plant that is controlled over a wireless network. 
Such a heterogeneity in the system structure often arises in various safety-critical applications which require different levels of system dynamics for handling various operational objectives.
One example of such applications is an NCS which consists of assembling manipulators and autonomous forklift in a manufacturing automation system \cite{de2006use}. 
This section describes and formally models the dynamics of each subsystem in the framework. 

\subsection{Markov Decision Process}
MDP models have been widely used to abstract the high level  dynamics of stochastic control systems. 
Formally stated, an MDP is a tuple 
\begin{equation}
\label{eq:mdp}
\mathcal{M}=(S, S_{\text{init}}, A, P, c_{M})
\end{equation}
where $S=\{1,2,\ldots, N\}$ is a finite set of states, $S_{\text{init}}$ is a finite set of initial states and $A=\{a_{1}, a_{2}, \ldots, a_{M_{a}}\}$ is a finite set of actions in which $M_a\in\mathbb{N}$ denotes the number of such actions. 
The matrix $P=\{p(i \vert j, a)\}_{i, j \in S, a \in A}$ is a transition matrix where each of its elements $p(i \vert j, a)$ denotes the transition probability from state $j$  to state $i$ with $(i,j)\in S$ under action $a\in A$. 
$c_{M}=\{c_{M}(s, a)\}_{s \in S, a \in A}$ is the set of costs where each of its elements $c_{M}(s, a) \in \mathbb{R}_{+}$ denotes the cost induced by $s\in S$ under $a\in A$.

As shown in Fig. \ref{fig: system-framework}, the proposed modeling framework uses the MDP  to describe stochastic discrete event systems that are managed to achieve certain high level tasks (e.g. transporting products between different workstations). 
In particular, these high level tasks are often required to be accomplished with certain levels of efficiency (e.g. with minimum time or energy). 
One main objective of this paper is to develop optimal control policies for such an MDP which ensure the overall system attains the specified level of efficiency. 

The locations of moving objects (e.g. cranes, forklift trucks and ground vehicles) in industrial workspace can significantly vary the channel characteristics of radio wave transmission and as a result causes the workspace to have different shadow fading levels  \cite{agrawal2014long, quevedo2013state,kashiwagi2010time}. 
In this regard, an MDP model  of a workspace can be constructed by partitioning the whole operational region into a finite number of small subregions with different fading levels. 
\subsection{Nonlinear Networked Control System}
\label{subsec:NNCS}
We now describe the model for the lower level tasks within the framework.
Consider a networked control system whose dynamics satisfy a nonlinear difference equation of the form
\begin{equation}
\label{eq: ncs}
x_{k+1}=f(x_k, u_k, w_k)
\end{equation}
where $x(\cdot): \mathbb{Z}_{\geq 0} \rightarrow \mathbb{R}^{n}$ is the vector of system states, $w(\cdot): \mathbb{Z}_{\geq 0} \rightarrow \mathbb{R}^{m}$ is a vector of essentially bounded external disturbances (i.e. $\exists M_w < \infty:|w|_{\mathcal{L}_{\infty}} \leq M_w$), and $u(\cdot): \mathbb{Z}_{\geq 0} \rightarrow \mathbb{R}^{l}$ is a vector of control input signals generated by a \emph{remote controller}. 
Given the signals $x(k), w(k), u(k)$ at time $k$, the function $f(\cdot): \mathbb{R}^{n}\times \mathbb{R}^{l}\times \mathbb{R}^{m}\rightarrow\mathbb{R}^n$  is a nonlinear map which determines the states $x_{k+1}$ of the system at the next time $k+1$.

As illustrated in Fig. \ref{fig: system-framework}, the states of the system are sampled at each discrete time $k$, encoded as a single data packet and then transmitted over a wireless fading channel. 
Let $\Omega_{p}=\{p_{i}\}_{i=1}^{M}$ be the set of transmission power levels that can be used to transmit the data packet. 
We assume that at each discrete time $k$, the data packet is transmitted using one communication power level $p_{k} \in \Omega_{p}$. 
Furthermore, given the inputs $x_{k}$ and $p_{k}$ to the fading channel, we assume that the data packet is received and decoded successfully at the end side of the channel with a  probability of $q_{k}$. 
We assume the probability $q_{k}$ is time-varying and varies as a function of the transmission power level $p_k$ and the MDP states $s_k$.

Now consider a sequence of independent random process $\{\gamma_{k}\}_{k=0}^{\infty}$ which may takes either $0$ or $1$ value at each time $k$. 
Let $\gamma_{k}$ be a random variable which characterizes the packet drop event at time $k$ of the form
\begin{equation}
\label{eq:gamma}
\gamma_{k}=\begin{cases}
1, &  \text{packet is successfully received and decoded,} \\
0, &     \text{packet is dropped.}
\end{cases}
\end{equation}
The probability measure of the random process $\{\gamma_{k}\}_{k=0}^{\infty}$ thus characterizes the impact that the fading channel  has on the successful reception of the packet. 
In this paper, a fading channel is called a \emph{state-dependent dropout channel} (SDDC) if it satisfies the condition in Definition \ref{def: SDDC} below.
\begin{definition}[SDDC]
	\label{def: SDDC}
	Consider a random process $\{\gamma_{k}\}_{k=0}^{\infty}$, an MDP $\mathcal{M}$ with a set of states $S$, and a set of transmission powers $\Omega_{p}$. 
	A fading channel is an SDDC if for all $i \in S$ and $p \in \Omega_{p}$ it has the following form of \emph{drop-out probability}.
	\begin{equation}
	\label{eq: SDDC}
	\mathbb{P}\{\gamma_{k}=0 \vert s_{k}=i, p_{k}=p\}=\theta(i, p), \;\;\forall k \in \mathbb{Z}_{\geq 0},
	\end{equation}
	where $\theta(i, p) \in [0, 1)$ is a continuous function of the MDP's state $i \in S$ which monotonically decreases with respect to the transmission power $p \in \Omega_{p}$.
\end{definition}

\begin{remark}
Note that the drop-out probability in \eqref{eq: SDDC} is closely related to the notion of \emph{outage probability} which is a widely used metric for characterizing the performance of the fading channels \cite{tse2005fundamentals}. 
Specifically, the outage probability is often used to measure the probability of receiving a data packet with a signal-to-noise ratio (SNR) below a prescribed threshold.
	In particular, let $N_{0}$ and $p_{k}$ denote the noise power and the transmission power at time $k$, respectively. 
	Let $h_{k}$ be a sequence of random variables which characterizes the variation of the channel gain at time $k$.
	In the case of factory automation systems, such a variation may occurs due to (i) multi-path propagation and (ii) shadow fading caused by moving obstacles (i.e. the states of the MDP).
	In this regard, we define $h_{k}=h^{p}_{k}h^{s}_{k}$ where $h^{p}_{k}$ and $h^{s}_{k}$ are  independent random variables, each of which denotes the variation induced by multi-path propagation and shadow fading, respectively.
	We further assume that $h^{p}_{k}$ is an independent, identically distributed random process with mean $\overline{h}^p$, whereas $h^{s}_{k}=\kappa(s_{k})$ is a state-dependent Markovian random process defined by a shadowing function $\kappa(\cdot): S \rightarrow \mathbb{R}_{+}$ which characterizes the impact that the MDP states have on the channel gain.
	From an industrial application standpoint, the MDP states $S$ in  (\ref{eq:mdp}) may represents the partition of the physical region of a factory/warehouse where the autonomous forklifts travel/operate.
	As such, one may define the shadowing function $\kappa$ by examining the fading level at each partition of the region. 
	In this regard, for a defined threshold level $\gamma_{0}$, the outage probability may be defined as
	\begin{align}
	\mathbb{P}\{\text{SNR} \leq \gamma_{0} &\,\vert\, s_{k}=i, p_{k}=p\} \nonumber\\
	&=\mathbb{P}\{(h^{p}_{k}h^{s}_{k})^{2}p_k \leq N_{0}\gamma_{0} \,\vert\, s_{k}=i, p_{k}=p\} \nonumber \\
	&=\mathbb{P}\left\{h^{p}_{k} \leq \frac{\sqrt{N_{0}\gamma_{0}}}{\sqrt{p}\kappa(s_k)} \,\bigg |\, s_{k}=i, p_{k}=p \right\} \nonumber \\
	&=\theta(i, p)
	\label{eq:drop-out}
	\end{align}
	Note that if the distribution $h^{p}_{k}$ of the channel is known in advance, one may explicitly define the expression of  $\theta(i, p)$. 
	For instance, if $h^{p}_{k}$ is known to have a Raleigh distribution with a scale parameter $\sigma=\overline{h}_{p}$, then $\theta(i, p)=1-e^{-\frac{N_{0}\gamma_{0}}{2p\kappa^{2}(i)\overline{h}_{p}^{2}}}$.
\end{remark}

Given the SDDC model \eqref{eq:gamma}-\eqref{eq:drop-out}, let $K(\cdot): \mathbb{R}^{n} \rightarrow \mathbb{R}^{l}$ be a controller for the NCS \eqref{eq: ncs}.
We assume the control signal $u_{k}$ in \eqref{eq: ncs} is constructed based on an the state estimate of the form
	\begin{align}
	\hat{x}_{k}&=\gamma_{k}x_{k}+(1-\gamma_{k})\hat{x}_{k-1}, \quad \hat{x}_{0}=0, \label{eq: estimator-z}\\
	u_{k}&=K(\hat{x}_{k}),
	\label{eq: controller-u}
	\end{align}
	where $\hat{x}_{k} \in \mathbb{R}^{n}$ denotes the estimate of the system state $x_{k}$ at time  $k$. 
	Note that \eqref{eq: estimator-z} is a zero-order hold estimator that holds the latest received information. 
	Now define an augmented state vector $z_{k}=[x_{k}; \,\hat{x}_{k-1}]$.
	Then the closed-loop system defined by \eqref{eq: ncs}, \eqref{eq: estimator-z} and \eqref{eq: controller-u} can be described as a randomly switched nonlinear system of the form
	\begin{align}
	\label{eq: randomly-switched}
	z_{k+1}=f_{\gamma_{k}}(z_{k}, w_{k}).
	\end{align}

Equation \eqref{eq: randomly-switched} basically suggests that the stability of the heterogeneous system model arises from \eqref{eq: ncs}, \eqref{eq: estimator-z} and \eqref{eq: controller-u} can be investigated using model \eqref{eq: randomly-switched}. 
One important challenge in such an analysis is the complexity of analyzing the cyber-physical coupling between the lower level plant model and the higher level MDP model in the presence of random switching signals $\{\gamma_{k}\}_{k=0}^{\infty}$ in \eqref{eq:gamma}-\eqref{eq: SDDC}. 
To ensure the overall system stability and efficiency, the transmission power and the control policies must thus be carefully designed.
%

\section{Stochastic Safety and Efficiency}
\label{sec: problem-formulation}
This section formally defines the notions of \emph{stochastic safety} and \emph{efficiency} which will be used in developing our proposed co-design policies framework for the heterogeneous systems in Fig. \ref{fig: system-framework}.
 Each of such definitions are detailed below.
	\begin{definition}[Stochastic Safety \cite{kozin1969survey}]
		\label{def: safety}
		Consider the NCS in \eqref{eq: ncs}, \eqref{eq: estimator-z}-\eqref{eq: controller-u} or \eqref{eq: randomly-switched}.
		let $\Omega_{s}=\{ x \in \mathbb{R}^{n} : |x| \leq r \}$ be a bounded safe region whose size is determined by a constant $r >0$. Define $\Omega_{r'}=\{ x \in \mathbb{R}^{n} : |x| <  r' \}$ to be a bounded set with $r' > 0$.  
		\begin{itemize}
			\item[\textbf{E1}] Let $w=0$ holds in \eqref{eq: ncs}, then \eqref{eq: ncs} is said to be \emph{asymptotically safe in expectation} (ASE) with respect to  $\Omega_{s}$ if there exist a class $\mathcal{KL}$ function $\beta(\cdot, \cdot)$ and a bounded set $\Omega_{r'}=\{ x \in \mathbb{R}^{n} : |x| <  r'(r) \}$ such that for all $x_{0} \in \Omega_{r'}$:
			\begin{align}
		 \mathbb{E}(|x_k|) \leq \beta(|x_0|, k) \leq r,   \quad \forall k \in \mathbb{Z}_{\geq 0}
			\end{align}
			and $\lim_{k \rightarrow \infty}\mathbb{E}(|x_k|)=0$.
			\item[\textbf{P1}] Let $w=0$ holds in \eqref{eq: ncs}. Then \eqref{eq: ncs} is said to be \emph{almost surely asymptotically safe} (ASAS) with respect to $\Omega_{s}$ if $\forall \epsilon \in (0, 1], k' > 0$,  there exists a class $\mathcal{KLL}$ function $\xi(\cdot, \cdot, \cdot)$ and a bounded set  $\Omega_{r'}=\{ x \in \mathbb{R}^{n} : |x| <  r'(\epsilon, r) \}$ such that $\forall x_{0} \in \Omega_{r'}$:
			\begin{align}
		\mathbb{P}\{\sup_{k \geq k'}|x_{k}| \geq r\} \leq \xi(|x_0|, k, r) \leq \epsilon
			\end{align}
			and ${\rm Pr}\{\lim_{k' \rightarrow \infty}\sup_{k \geq k'}|x_{k}| \geq r\} = 0$.
			\item[\textbf{P2}] The system in \eqref{eq: ncs} with ultimately bounded disturbance $|w|_{\mathcal{L}_{\infty}} \leq M_w$ is said to be \emph{practically safe in probability} (PSP) if  for any $\epsilon > 0$ and $ \delta > 0$, there always exists $\rho_{\epsilon}(\Delta, M_w) \in (0, 1)$ such that for all $x_{0} \in \Omega_{s}$,
			\begin{align}
			\label{ineq: practical-safe-in-probability}
			\lim_{k \rightarrow \infty}\mathbb{P}\{|x_{k}| \geq \epsilon+\Delta\} \leq \rho_{\epsilon}(\Delta, M_w).
			\end{align}
		\end{itemize}
\end{definition}

\begin{remark}
	Note that the stochastic safety is defined in terms of the infinite norm $|x|$ of the system states.
	It can be shown that the ASE property implies \emph{mean square safety} (MSE) that characterizes the variance of the system states (i.e., $\mathbb{E}(\|x_k\|_2) \leq r, \forall k \in \mathbb{Z}_{\geq 0}$) by realizing that $\|x\|_{2} \leq \sqrt{n}|x|$. 
	Notice also that \textbf{E1} and \textbf{P1} in Definition \ref{def: safety} suggest that ASAS is stronger notion than ASE. 
	Section \ref{sec: safety} particularly shows that sufficient conditions to assure \textbf{P1} are indeed stricter than \textbf{E1}.  
\end{remark}

On the one hand, the safety of system \eqref{eq: ncs} can be achieved by enforcing sufficiently small outage probability. 
As suggested in the SDDC model \eqref{eq: SDDC}, this can be done by adjusting both the control policies of the MDP model and the transmission power $p$ of the communication system. 
On the other hand, the efficiency of the system may be maintained by designing optimal control and power policies which optimize some predefined long-run performance criteria. 

The control policies for both the MDP and transmission power models are defined as infinite sequences $\mu^{m}=(\mu_{1}^{m}, \mu_{2}^{m}, \ldots)$ and $\mu^{p}=(\mu^{p}_{1}, \mu^{p}_{2}, \ldots)$, respectively. 
Let $\mu_{k}=(\mu^{m}_{k}, \mu^{p}_{k})$ denotes the randomized decision rule  at time $k$. 
The design of such a decision rule may depends on the  information observed along the time. 
	Specifically,  let $h_{k}=(s_{1}, a_{1}, p_{1}, \ldots, s_{k-1}, a_{k-1}, p_{k-1}, s_{k})$ denotes the history of the states, actions and power up to time $k$.
	Then the decision variable $\mu_{k}$  at time $k$ is defined as the conditional probability distribution over the sets of  action $A$ and transmission power $\Omega_{p}$, given the history information $h_{k}=(h_{k}^{m}, h_{k}^{p})$ up to time $k$.
	More formally, for all $a \in A$ and $p' \in \Omega_{p}$, we define
	\begin{equation*}
	\mu_{k}^{m}(a)={\rm Pr}\{a \vert h_{k}\}\qquad \text{and}\qquad \mu_{k}^{p}(p')={\rm Pr}\{p' \vert h_{k}\}.
	\end{equation*}
	In particular, each of such decision variable or policy $\mu_{\infty}^{\ell}$ for $\ell \in \{p,m\}$ is said to be \emph{stationary} if the corresponding decision making variable is time homogeneous and depends only on the states, i.e., $\mu_{\infty}^{p}={\rm Pr}\{p \vert s\}, \mu_{\infty}^{m}={\rm Pr}\{a \vert s\}, \forall s \in S$. 
	We say that the heterogeneous system in Fig. \ref{fig: system-framework} is \emph{stochastically efficient} under such stationary policy if Definition \ref{def: stochastic-efficiency} below holds.

\begin{definition}[Stochastic Efficiency]
		\label{def: stochastic-efficiency}
		Given the set of transmission powers $\Omega_{p}=\{p_{\ell}\}_{\ell=1}^{M}$ and the set of costs $c_{M}=\{c_{M}(s, a)\}_{s \in S, a \in A}$, let $c_{p}=\{c_{p}(p)\}_{p \in \Omega_{p}}$ be the set of costs induced by $\Omega_{p}$. 
		We say that the heterogeneous system in \eqref{eq: ncs}, \eqref{eq: estimator-z}-\eqref{eq: controller-u} or \eqref{eq: randomly-switched} is \emph{stochastically  efficient} with respect to $c_{M}$ and $c_{p}$ if it minimizes an infinite time average cost $J(\cdot)$ below.
		\begin{equation}
		\label{eq: opt-efficiency}
		J(s_{0}, \mu^{m}, \mu^{p})=\lim_{T \rightarrow \infty}\frac{1}{T} \mathbb{E}_{s_{0}}^{\mu^{m}, \mu^{p}} \left\{\sum_{i=0}^{T} c_{M}(s_{k}, a_{k})+\lambda c_{p}(p_{k}) \right\}
		\end{equation}
		where $\lambda > 0$ is a constant and $\mu^{m}=\{\mu^{m}_{0}, \mu^{m}_{1}, \ldots\}$ and $\mu^{p}=\{\mu^{p}_{0}, \mu^{p}_{1}, \ldots\}$ denote control and power policies, respectively. 
\end{definition}

The \emph{cyber-physical coupling} between the NCS \eqref{eq: ncs} and MDP \eqref{eq:mdp} suggests the need to solve a co-design problem to obtain control policies that achieve both the safety and efficiency of the heterogeneous system \eqref{eq: randomly-switched}.
	Sections \ref{sec: safety} and \ref{subsec: co-design} propose our solution framework to develop such a co-design problem.
	Our framework essentially formulates a constrained cooperative game problem that seeks to find equilibrium points consisting the optimal control and power strategies that will guarantee the overall system's safety and efficiency.
\section{Standing Assumptions}
\label{sec:assumptions}
This section states two assumptions which underlie the  main results in Sections \ref{sec: safety} and \ref{subsec: co-design}. 
The first one (Assumption \ref{assumption: multiple-lyapunov}) assumes the existence of a multiple Lyapunov functions (MLF) for the randomly-switched system model in \eqref{eq: randomly-switched}. 
The second one (Assumption \ref{assumption: unichain}) assumes the MDP \eqref{eq:mdp} satisfies the so-called \emph{unichain} structure \cite{puterman1994markov, anthonisse1977exponential} for any given control policy $\mu^{m}$. 
Both of these assumptions are formally stated below.

\begin{assumption}[Existence of MLF]
	\label{assumption: multiple-lyapunov}
	Consider the randomly switched system description in \eqref{eq: randomly-switched}. 
	Let $\{V_{i}\}_{i \in \{0, 1\}}$ be a family of Lyapunov functions. 
	For system \eqref{eq: randomly-switched}, there exist class $\mathcal{K}_{\infty}$ functions  $\alpha_{j}(\cdot)$, class $\mathcal{K}$ functions $\chi(\cdot)$, and constants $\varrho, \lambda_{i} > 0$ for  $i=(0, 1)$ and $j=(1, 2)$ such that
	\begin{itemize}
		\item[(V1)] $\alpha_{1}(|x|) \leq V_{i}(x) \leq \alpha_{2}(|x|), \;\; \text{for all }\; i \in \{0, 1\}$	
		\item[(V2)] $V_{i}(f_{i}(x, w)) \leq \lambda_{i} V_{i}(x)+\chi(|w|),  \;\; \text{for all }\; i \in \{0, 1\}$
		\item[(V3)] $V_{i}(x) \leq \varrho V_{j}(x),\;\;\text{ for all }\; i, j \in \{0, 1\} $ with $i \neq j$
		\item[(V4)] $\varrho \min_{i}\,\lambda_i < 1$.
	\end{itemize}
\end{assumption}
\begin{remark}
In Assumption \ref{assumption: multiple-lyapunov}, condition V1  essentially requires the MLF to be radially unbounded whereas condition V2 necessitates the growth of the $i$th Lyapunov function $V_i(\cdot)$ of the MLF along the vector field $f_{i}(x, w)$ of the $i$th subsystem to be bounded from above by the combination of a linear function of $V_{i}(x)$ and a function of disturbance $|w|$. 
When these conditions are satisfied with $\lambda_{i} < 1$, the $i$th subsystem $f_{i}(x, w)$ of  \eqref{eq: randomly-switched} is said to be \emph{discrete input-to-state stable}   \cite{jiang2001input}. 
Note that  condition V3  is a commonly used  assumption when proving the stability of  switched systems (cf. e.g. \cite{liberzon2012switching}) even though it restricts the class of Lyapunov functions that may be used as the MLF. 
For instance, quadratic Lyapunov functions in switched linear systems stability analysis certainly satisfy V3. 
Finally, condition V4 requires that at least one subsystem of \eqref{eq: randomly-switched} is sufficiently stable (can be enforced by a controller which guarantees a sufficiently small decay rate $\lambda_{i} < 1$  when the communication is perfect  without packet dropout  ($\gamma_k=1$)). 
\end{remark}

\begin{assumption}[Unichain MDP \cite{puterman1994markov, anthonisse1977exponential}]
	\label{assumption: unichain}
	Let $P_{k}:=P(\mu_{k}^{m})$ be the transition probability matrix~(stochastic matrices) of the MDP under the policy $\mu_{k}^{m}$ at time $k$. For a given $n > 0$, $\forall 1 \leq k \leq n$, 
the products of the stochastic matrices, $\Pi_{i=1}^{k}P_{i}$ induced by any selected control policy $\mu^{m}=(\mu_{1}^{m}, \mu_{2}^{m}, \ldots, \mu_{k}^{m})$, is aperiodic and has a single ergodic class.
\end{assumption}

Under the unichain condition in Assumption \ref{assumption: unichain}, it was shown in  \cite[Theorem 1]{anthonisse1977exponential} that the Markov chain\footnote{Possibly time-inhomogeneous Markov chain} induced by the control policy $\mu^{m}$ exponentially converges to a stationary distribution. 
This fact is recalled below for completeness. 
\begin{theorem}[\cite{anthonisse1977exponential}]
	\label{thm:old}
	Suppose Assumption \ref{assumption: unichain} holds and let there be an integer $v\geq 1$ and a real number $0 \leq \alpha < 1$.
	Then for any sequence $\{P_i, i \geq 1\}$ of stochastic matrices, there is a probability distribution $\{\pi_{i}, 1 \leq i \leq N\}$ such that for all $i$ 
	\begin{align}
		\label{ineq-convergence}
	|(P_{n}\cdots P_{1})_{ij}-\pi_{i}| \leq \alpha^{[n / v]}, \quad \forall n \geq 1, j=1,2,\ldots, N,
	\end{align}
	where $[x]$ is the largest integer that is less than or equal to $x$.
\end{theorem}

\section{Stochastic Safety}
\label{sec: safety}
This section presents the result of this paper on stochastic safety of the heterogenous system in Fig. \ref{fig: system-framework}. 
Firstly, when external disturbance $w$ is absent in NCS \eqref{eq: ncs}, sufficient conditions which ensure ASE~(i.e. \textbf{E1} in Definition \ref{def: safety}) and ASAS~(i.e. \textbf{P1} in Definition \ref{def: safety}) are derived in Theorems \ref{thm: safety-expectation} and \ref{thm-asas}, respectively.
When bounded disturbance is present, sufficient conditions which ensure a weaker notion of PSP~(i.e. \textbf{P2} in Definition \ref{def: safety}) is then presented in Theorem \ref{thm: practical-safety}.

As stated in Section \ref{sec: problem-formulation}, the safety problem of the NCS \eqref{eq: ncs} can be examined using randomly switched systems framework \cite{chatterjee2011stabilizing}. 
Under Assumption \ref{assumption: multiple-lyapunov}, the safety of the switched systems in \eqref{eq: randomly-switched} is closely related to the probability measure of the switching signal $\{\gamma_{k}\}$. 
As shown in \eqref{eq: SDDC}, the distribution of the random processes $\{\gamma_{k}\}$ is governed by the states of both MDP model $S$ and transmission power set $\Omega_{p}$. 
In order to ensure the safety of the randomly switched system \eqref{eq: randomly-switched}, both the MDP and transmission power must be coordinated to ensure an appropriate behavior of the switching signal $\{\gamma_{k}\}$. 

Let $\overline{s}_k=(s_k, p_k)$ be a random variable representing the joint-state of MDP system and transmission power at time $k$.
Proposition \ref{proposition: transition-probability} shows that the corresponding random process $\{\overline{s}_{k}\}$ that is defined on the joint state space $S \times \Omega_{p}$ is a Markov process whose transition probability matrix is a function of the transmission power policy $\mu^{p}$, transition matrix $P$ and control policy $\mu^{m}$ of the MDP.   
\begin{proposition}
\label{proposition: transition-probability}
Consider the MDP $\mathcal{M}$ and a finite set of transmission power $\Omega_{p}=\{p_{i}\}_{i=1}^{M}$.
For a joint policy $\mu=(\mu^{p}, \mu^{m})$ on the joint state space $\overline{S}:=S \times \Omega_{p}$, the random process $\{\overline{s}_{k}\}$ with $\overline{s}_{k}=(s_{k}, p_{k})$ is a Markov process.
The transition probability matrix $\overline{P}_{k}:=\overline{P}(\mu_{k})$ of $\{\overline{s}_{k}\}$ at time $k$ is defined over $\overline{S}$ and $\mu_k$ and takes the form
\begin{align}
\label{eq: P_k}
\overline{P}_{k}(\overline{s}, \overline{s}')&=\mathbb{P}\{\overline{s}_{k+1}=\overline{s}' \vert \overline{s}_{k}=\overline{s}\} \nonumber \\
&=\mu_{k}^{p}(p')\sum_{a \in A(s)}p(s' \vert s, a) \mu_{k}^{m}(a), \quad \forall \overline{s}', \overline{s} \in \overline{S}.
\end{align}
where $\overline{s}=(s, p), \overline{s}'=(s', p') \in \overline{S}$ whereas $\mu_{k}^{p}(p')=\mathbb{P}\{p' \vert \overline{s}_{k}=s\}$ and $\mu^{m}_{k}(a)=\mathbb{P}\{a \vert s_{k}=s \}$.
\end{proposition}
\begin{IEEEproof}
Let $\Omega^{k}:=\underbrace{S \times \Omega_{p} \times  A \times \cdots \times S \times \Omega_p \times A}_{=k-1} \times S \times \Omega_{p}$ denotes the sample space up to time $k$. 
The random process $\{(s_k, p_k)\}$ is defined over $\Omega^{k}$ such that for any sample path $h_{k}=(s_{1}, p_1, a_1, \ldots, s_{k}, p_{k}) \in \Omega^k$, one has
\begin{align}
\label{eq:add1}
\mathbb{P}\{s_{k+1}=&s', p_{k+1}=p' \vert h_{k}\} \nonumber \\ 
&=\mathbb{P}\{s_{k+1}=s' \vert h_{k}\}\mathbb{P}\{p_{k+1}=p' \vert s_{k+1}=s', h_{k}\} \nonumber \\
&=\mathbb{P}\{s_{k+1}=s' \vert s_{k}, p_k \}\mathbb{P}\{p_{k+1}=p' \vert s_{k+1}=s', s_k, p_k\} \nonumber\\
&=\mathbb{P}\{s_{k+1}=s', p_{k+1}=p' \vert s_k, p_k\}
\end{align}
The first and third equalities in \eqref{eq:add1} are based on the chain rule of probability. 
The second equality holds because the random process $\{s_k\}$ is a Markov process with a Markovian transmission power policy $\mu^{p}$.
As such, the random process $\{\overline{s}_{k}\}$ is Markovian by definition. 
Now for any joint-states $\overline{s}'=(s', p'), \overline{s}=(s, p) \in S \times \Omega_p$, the transition probability from $\overline{s}$ to $\overline{s}'$ under the Markovian policy $\mu_k=(\mu^{p}_{k}, \mu^{m}_{k})$ at time $k$ can then be constructed and is given by
\begin{align*}
\mathbb{P}\{\overline{s}_{k+1}=\overline{s}' &\vert \overline{s}_{k}=\overline{s}\} \\
&=\mathbb{P}\{p_{k+1}=p' \vert \overline{s}_{k}=\overline{s} \}\times\\
&\qquad \quad  \mathbb{P}\sum_{a \in A(s)}\{s_{k+1}=s' \vert s_{k}=s, a\} \mathbb{P}\{a \vert s_{k}=s\}\\
&=\mu_{k}^{p}(p')\sum_{a \in A(s)}p(s' \vert s, a) \mu_{k}^{m}(a)
\end{align*}
as claimed in \eqref{eq: P_k}. The proof is thus completed. 
\end{IEEEproof}

Under Assumption \ref{assumption: unichain}, Theorem \ref{thm: safety-expectation} below sates sufficient conditions that guarantee the NCS \eqref{eq: randomly-switched} is  ASAS in the absence of external disturbance~(i.e., $w=0$) or PSP when an ultimately bounded disturbance is present~(i.e., $0 < |w|_{\mathcal{L}_{\infty}} < \infty$). 

\begin{theorem}
	\label{thm: safety-expectation}
Consider the NCS (or randomly nonlinear switched system) in \eqref{eq: randomly-switched} without external disturbance ($w=0$) and the SDDC model in \eqref{eq: SDDC}.
Suppose Assumptions \ref{assumption: multiple-lyapunov} and \ref{assumption: unichain} hold.
For a selected joint-policy $\mu=\{\mu_{k}\}_{k=1}^{\infty}$, let $\{\overline{P}_{k}, k \geq 1\}$ be the transition matrix of the Markov chain that is induced by the joint policy $\{\mu_{k}, k \geq 1\}$.
Let $\{\overline{\pi}_{i}, 1 \leq i \leq NM\}$ denotes the stationary distribution for the steady states of the Markov chain.
Then the NCS in \eqref{eq: randomly-switched} is ASE with respect to the origin if the condition below holds.
\begin{align}
\label{ineq:sufficient-condition}
\theta^{T}\overline{\pi} < \frac{1-\lambda_{1}\varrho}{\varrho(\lambda_{0}-\lambda_{1})}
\end{align}
where $\theta=[\theta(s, p)]_{NM \times 1}$ and $\overline{\pi}=[\mathbb{P}\{s, p\}]_{1 \times NM}$. Also, if $\exists \eta \in [0, 1)$ such that the condition $\theta^{T}\overline{\pi} < (\eta-\lambda_{1}\varrho)\big/[\varrho(\lambda_{0}-\lambda_{1})]$ holds, then the NCS is \emph{exponentially safe in expectation}, i.e., $\exists \kappa > 0$ and a function $\alpha\in\mathcal{K}_{\infty}$ such that $\mathbb{E}(|x_k|) \leq \kappa \eta^{k}\alpha(|x_0|)$.
\end{theorem}
\begin{IEEEproof}
	See Appendix \ref{appendix:proof}.
\end{IEEEproof}
\begin{remark}
Condition \eqref{ineq:sufficient-condition} essentially implies that the NCS \eqref{eq: randomly-switched} is guaranteed to be ASE only if the state-dependent channel (of packet loss distribution $\theta$) is sufficiently `good' on average.
\end{remark}
\begin{remark}
By Markov inequality, it can be shown that \eqref{ineq:sufficient-condition} ensures the NCS \eqref{eq: randomly-switched}  is \emph{stochastically safe with probability one} \cite{kozin1969survey}.
This basically means that for any initial state $x_0\in\Omega_{s}$, there always exists a constant $\epsilon > 0$ such that $\lim_{k \rightarrow \infty}{\rm Pr}\{|x_k| \geq \epsilon \} = 0$. 
The  notion of safety in probability is weaker than that of ASAS (cf. Definition \ref{def: safety}) in that the former only requires the probability of the system states leaving $\Omega_{s}$ goes to zero in the time limit.
 In contrast,  the notion of ASAS requires that the probability of almost all of the system's sample paths with $x_0\in\Omega_{s}$  exiting the safe region is arbitrarily small for any infinite time duration $[k', \infty)$ and goes to zero as time $k'$ goes to infinity.  
 It is clear that a stronger condition than \eqref{ineq:sufficient-condition} is needed to ensure almost sure safety.
\end{remark}

Theorem \ref{thm-asas} below gives a sufficient condition to ensure the NCS \eqref{eq: ncs} is ASAS in the absence of external disturbance. 
\begin{theorem}
\label{thm-asas}
Consider the NCS \eqref{eq: ncs} without external disturbance ($w=0$) and the SDDC model in \eqref{eq: SDDC}.
Suppose Assumption \ref{assumption: multiple-lyapunov} holds.
Then for a selected joint policy $\mu=\{\mu_{k}, k \geq 1 \}$,  the NCS \eqref{eq: ncs} is ASAS if $\forall k \in \mathbb{Z}_{\geq 0}$ and  $\forall (s, p) \in S \times \Omega_p$:
\begin{align}
\label{ineq: sufficient-safety}
\sum_{\substack{s' \in S \\  p' \in \Omega_p}}\theta(s', p')\mu_{k}^{p}(p')\sum_{a \in A(s)} p(s' \vert s, a)\mu_{k}^{m}(a)< \frac{1-\lambda_{1}\varrho}{\varrho(\lambda_{0}-\lambda_{1})}
\end{align}
where $\mu_{k}^{p}(p')=\mathbb{P}\{p' \vert s_{k}=s\}$ and $\mu_{k}^{m}(a)=\mathbb{P}\{a \vert s_{k}=s\}$. 
Moreover, if $\exists \eta \in [0, 1)$ such that the condition in \eqref{ineq: sufficient-safety} holds with the threshold $(\eta-\lambda_{1}\varrho)\big/[\varrho(\lambda_{0}-\lambda_{1})]$, then the NCS is \emph{almost surely exponential safe}, i.e., $\exists \kappa(r) > 0$ and a class $\mathcal{K}_{\infty}$ function $\alpha(\cdot)$ such that $\mathbb{P}\{\sup_{k \geq k'}|x_{k}| \geq r\} \leq \kappa(r) \eta^{k} \alpha(|x_0|)$.
\end{theorem}
\begin{IEEEproof}
	See Appendix \ref{appendix:proof}.
\end{IEEEproof}

Corollary \ref{coll:add1} below states an inference based on the derived conditions in \eqref{ineq:sufficient-condition} and \eqref{ineq: sufficient-safety}.
\begin{corollary}
\label{coll:add1}
	For the NCS \eqref{eq: ncs} with no disturbance ($w=0$) and the SDDC model in \eqref{eq: SDDC}, the condition to ensure ASAS is stronger than that of ASE in the sense that \eqref{ineq: sufficient-safety} implies \eqref{ineq:sufficient-condition}.
\end{corollary}
\begin{IEEEproof}
By \eqref{eq: P_k}, the  inequality in \eqref{ineq: sufficient-safety} is equivalent to
\begin{align*}
 \sum_{\overline{s}' \in S \times \Omega_{p}}\theta(\overline{s}')\overline{P}_{k}(\overline{s}, \overline{s}') < \frac{1-\lambda_{1}\varrho}{\varrho(\lambda_{0}-\lambda_{1})}, \; \forall \overline{s} \in S \times \Omega_{p}.
\end{align*}
and thus may be rewritten in a vectorial form below
\begin{align*}
\theta^{T}\overline{P}_{k} < \frac{1-\lambda_{1}\varrho}{\varrho(\lambda_{0}-\lambda_{1})} \bm{1}, \; \forall k \in \mathbb{Z}_{+}
\end{align*}
where $\bm{1}=[1, 1, \ldots, 1]_{1 \times NM}$ is an all-one vector and $\overline{P}_{k}$ is the transition matrix for the joint-states in \eqref{eq: P_k}. 
By Theorem \ref{thm:old} and Assumption \ref{assumption: unichain}, we know that there exists a stationary distribution $\overline{\pi}$ over the joint state set $S \times \Omega_{p}$ such that $\prod_{k=1}^{\infty}\overline{P}_{k}\overline{\pi}=\overline{\pi}$. 
For a given $\ell \in \mathbb{Z}_{+}$, suppose for any $k \leq \ell$, the condition \eqref{ineq: sufficient-safety} holds such that
\begin{align}
\label{ineq: multiple-P}
\theta^{T}\prod_{k=1}^{\ell} \overline{P}_{k}< \frac{1-\lambda_{1}\varrho}{\varrho(\lambda_{0}-\lambda_{1})} \bm{1}, \; \forall \ell \in \mathbb{Z}_{+}.
\end{align}
Note that  \eqref{ineq: multiple-P} holds because $\{\overline{P}_{k}, \forall k\}$ are left stochastic matrices with each column summing up to $1$.  
Since $\ell \in \mathbb{Z}_{+}$is arbitrarily selected, let $\ell \rightarrow \infty$.
Then the condition in \eqref{ineq:sufficient-condition} can be recovered by right multiplying both sides of the inequality \eqref{ineq: multiple-P} with the stationary probability vector $\overline{\pi}$, i.e., 
\begin{align*}
\theta^{T}\underbrace{\prod_{k=1}^{\infty} \overline{P}_{k}\overline{\pi}}_{=\overline{\pi}} < \frac{1-\lambda_{1}\varrho}{\varrho(\lambda_{0}-\lambda_{1})} \underbrace{\bm{1} \overline{\pi}}_{=1} \Leftrightarrow \theta^{T}\overline{\pi} < \frac{1-\lambda_{1}\varrho}{\varrho(\lambda_{0}-\lambda_{1})}.
\end{align*}
Thus, the above inequality shows that condition in \eqref{ineq: sufficient-safety} implies the condition in \eqref{ineq:sufficient-condition}. The proof is therefore completed. 
\end{IEEEproof}

Next, Theorem \ref{thm: practical-safety} below derives a sufficient condition that guarantees the NCS \eqref{eq: ncs} is PSP.

\begin{theorem}
\label{thm: practical-safety}
Consider the MDP $\mathcal{M}$, a transmission power set $\Omega_p$ and the SDDC model in \eqref{eq: SDDC}. Suppose that Assumption \ref{assumption: multiple-lyapunov} holds.
Then the NCS in \eqref{eq: ncs}, \eqref{eq: estimator-z} and \eqref{eq: controller-u} with an ultimately bounded disturbance $w$ (i.e.,  $\exists 0 < M_w < \infty: \|w\|_{\mathcal{L}_{\infty}} \leq M_w$) is PSP if, for a given $\eta \in (\lambda_{1}\varrho, 1)$, there exists a joint policy $\mu=\{\mu_{k}, k \geq 1\}$ such that for all $ (s, p) \in S \times \Omega_p$
\begin{align}
\label{ineq: practical-stable-in-probability}
\sum_{\substack{s' \in S\\  p' \in \Omega_p}}\theta(s', p')\mu_{k}^{p}(p')\sum_{a \in A(s)} p(s' \vert s, a)\mu_{k}^{m}(a) \leq \frac{\eta-\lambda_{1}\varrho}{\varrho(\lambda_{0}-\lambda_{1})},
\end{align}
where $\mu_{k}^{p}(p')=\mathbb{P}\{p' \vert s_{k}=s\}$ and $\mu_{k}^{m}(a)=\mathbb{P}\{a \vert s_{k}=s\}$. 
In particular, let $\Omega_{s}=\{x \in \mathbb{R}^{n} \vert |x| \leq \Delta \}$ denotes the target region towards which the system state $x$ is driven.
Then the probability of exiting $\Omega_{s}$ as defined in \eqref{ineq: practical-safe-in-probability} satisfies
\begin{align}
\rho_{\epsilon}(\Delta, M_w)=\frac{\chi(M_w)}{(1-\eta)\alpha_{1}(\Delta+\epsilon)},
\label{eq: exiting-probability}
\end{align}
where $\chi(\cdot)$ and $\alpha_{1}(\cdot)$, respectively, are class $\mathcal{K}$ and $\mathcal{K}_{\infty}$ functions defined in Assumption \ref{assumption: multiple-lyapunov}.
\end{theorem}
\begin{IEEEproof}
See Appendix \ref{appendix:proof}.
\end{IEEEproof}
\begin{remark}
Note that the variable $\eta \in (\lambda_{1}\varrho, 1)$ in Theorem \ref{thm: practical-safety} is the convergence rate which characterizes how fast the expected value of the system state $x$ moves toward the target region $\Omega_{s}$. 
It is clear from \eqref{eq: exiting-probability} that the likelihood of the system state leaving the target region $\Omega_{s}$ is a monotonically decreasing function of the convergence rate $\eta$ and the disturbance magnitude $M_w$. 
Moreover, the exiting probability $\rho_{\epsilon}(\Delta, M_w)$ increases when the size of the target region reduces. 
\end{remark}
\begin{remark}
The safety condition in \eqref{ineq: practical-stable-in-probability} and exiting probability in \eqref{eq: exiting-probability} suggest that there exists a trade off between the system performance ($\eta, \rho_{\epsilon}(\Delta, M_w)$) and the joint control-communication policies ($\mu$) that are used to achieve that performance. 
In particular, inequality \eqref{ineq: sufficient-safety} implies that a higher convergence rate (or low exiting probability) leads to a smaller set of joint policies $\mu$ that are feasible to optimize the overall system costs defined in Definition \ref{def: stochastic-efficiency}.
\end{remark}
\section{Co-design of Safety and Efficiency: A Two-Player Constrained Cooperative Game}
\label{subsec: co-design}
This section considers a joint design framework to achieve both safety and efficiency of the heterogeneous NCS \eqref{eq: ncs}. 
In particular, it is shown that system efficiency can be assured if the infinite average cost in \eqref{eq: opt-efficiency} is minimized under the safety constraints \eqref{ineq: sufficient-safety} or \eqref{ineq: practical-stable-in-probability}. 
Based on the sufficient conditions derived in Section \ref{sec: safety}, it is shown  that the co-design problem of safety and efficiency for the whole industrial NCS can be formulated as a two-player constrained cooperative game.
In the formulated game, the previously derived sufficient conditions are posed as constraints on the strategy space that each player has to obey. 
Furthermore, Sections  \ref{subsec: lp} and \ref{subsec: polynomial-opt} also present efficient algorithms  to solve such a game.


Formally, the co-design of safety and efficiency can be formulated as the following constrained optimization problem.
\begin{problem}[Two-player Constrained Cooperative Game]
	\label{problem: constrained-cooperative-game}
Consider the notions of \emph{stochastic safety}  and \emph{stochastic efficiency}  in Definitions \ref{def: safety} and \ref{def: stochastic-efficiency}, respectively.
The solution to the co-design problem of the NCS \eqref{eq: ncs} is given by the joint policy $\mu=(\mu^{m}, \mu^{p})$ which solves the optimization below.
\begin{equation}
\label{eq:game}
\begin{aligned}
& \underset{\mu}{\min}
& & J(s_{0}, \mu^{m}, \mu^{p}) \; \text{in} \; \eqref{eq: opt-efficiency} &\\
& \text{s. t.}
& & \eqref{ineq: sufficient-safety} \; \text{or} \; \eqref{ineq: practical-stable-in-probability}&\qquad\text{(safety conditions) }
\end{aligned}
\end{equation}
\end{problem}

The main challenge in solving Problem \ref{problem: constrained-cooperative-game} lies in the structure of the safety conditions which impose constraints on the decision space of the players~(i.e. transmission power and MDP controller) in the cooperative game. 
Under the derived ASE condition in \eqref{ineq:sufficient-condition}, it is shown in Section \ref{subsec: lp} that the equilibrium of the game in Problem \eqref{problem: constrained-cooperative-game}~(i.e. optimal power and control policies) can be obtained by linear programming. 
When instead the stronger ASAS condition \eqref{ineq: sufficient-safety} is enforced as the constraint,  Section \ref{subsec: polynomial-opt} shows that the optimal solution of Problem \eqref{problem: constrained-cooperative-game} can be obtained using quadratic programming.  

\subsection{Co-design of stochastic efficiency and ASE: A linear programming approach} 
\label{subsec: lp}
Under Assumption \ref{assumption: unichain} and the  condition in \eqref{ineq:sufficient-condition}, Proposition \ref{proposition: stationary-policy} shows that the search for optimal solutions to Problem \ref{problem: constrained-cooperative-game} may sufficiently be done over the stationary polices.
\begin{proposition}
	\label{proposition: stationary-policy}
Let $U^{H}$ and $U^{S}$ denote, respectively the history-dependent policy and the stationary policy spaces for the MDP system $\mathcal{M}$ and the transmission power $\Omega_p$. 
Suppose Assumption \ref{assumption: unichain} holds.
Let $\mu^{*} \in U^{H}$ be an optimal solution to Problem \ref{problem: constrained-cooperative-game} under the safety constraint \eqref{ineq:sufficient-condition} and let $J^{*}(s_0, \mu^{*})$ be the corresponding optimal cost.
Then there always exists a stationary policy $\mu^{*}_{\infty} \in U^{S}$ such that $J^{*}(s_0, \mu^{*})=J^{*}(s_0, \mu^{*}_{\infty})$.
\end{proposition}
\begin{IEEEproof}
Note that Problem \ref{problem: constrained-cooperative-game} with a constraint as defined in \eqref{ineq:sufficient-condition} is equivalent to the constrained MDP problem with expected average costs (cf. \cite[Theorem 4.1]{altman1999constrained} for a proof of the completeness of stationary policies).  
As such, the proof of Proposition \ref{proposition: stationary-policy} is similar to that in \cite[Theorem 4.1]{altman1999constrained} and thus will only be sketched here. 
First, by Proposition \ref{proposition: transition-probability}, we know that the the joint state $\overline{s}_k=(s_k, p_k)$, $s_k \in S, p_k \in \Omega_p$ is an MDP with a transition probability as in \eqref{eq: P_k}.
Secondly, note that the ASE condition in \eqref{ineq:sufficient-condition} is a linear constraint which only depends on the stationary distribution of the joint state $S \times \Omega_{p}$. 
Thirdly, by the results of \cite[Theorem 4.1]{altman1999constrained}, it can be concluded that for any history-dependent policy which assures an optimal cost, there always exists a corresponding stationary policy that attains the same optimal cost and also satisfies the linear constraint \eqref{ineq:sufficient-condition}. The sketch of proof is completed. 
\end{IEEEproof}

Based on Proposition \ref{proposition: transition-probability}, a linear programming (LP) problem for computing the optimal solutions to Problem \ref{problem: constrained-cooperative-game} under the safety constraint \eqref{ineq:sufficient-condition} is formulated in \eqref{opt: lp} below. 
	\begin{align}
& \underset{\substack{\{X_{1}(s, a)\} \\ \{X_{2}(s, p)\}}}{\min}
& &\sum_{\substack{s \in S, a \in A(s) \\ p \in \Omega_{p}}} X_{1}(s, a)c_{M}(s, a)+\lambda X_{2}(s, p) c_{p}(s, p) \nonumber \\
& \text{s.t.}
& & \begin{cases}
& \sum\limits_{a \in A(s')} X_{1}(s', a)\\ 
&-\sum\limits_{\substack{s \in S \\ a \in A(s)}} p(s' | s, a)X_{1}(s, a)=0, \forall s' \in S \\
& \sum\limits_{\substack{s \in S \\ a \in A(s)}} X_{1}(s, a)=1,\sum\limits_{\substack{s \in S \\ p \in \Omega_{p}}}X_{2}(s, p)=1\\
& \sum\limits_{\substack{s \in S\\ p \in \Omega_{p}}} \theta(s, p)X_{2}(s, p) < \frac{1-\lambda_{1} \rho}{\rho(\lambda_{0}-\lambda_{1})}.
\end{cases}
\label{opt: lp}
\end{align}	
where $\{X_{1}(s, a)\}_{s\in S, a \in A(s)}$ and $\{X_{2}(s, p)\}_{s \in S, p \in \Omega_{p}}$ are the decision variables of the program. 
As discussed in \cite{altman1999constrained}, the decision variables $\{X_{1}(s, a)\}_{s\in S, a \in A(s)}$ and $\{X_{2}(s, p)\}_{s \in S, p \in \Omega_{p}}$ represent the stationary probability distributions over the state-action set $S \times A$ and the state-power set $S \times \Omega_{p}$, respectively, due to the stationary policy result in Proposition \ref{proposition: stationary-policy}. 
Theorem \ref{thm: lp} below shows that the solutions of the LP problem \eqref{opt: lp} define the optimal control and communication policies.
\begin{theorem}
	\label{thm: lp}
Consider the co-design Problem \ref{problem: constrained-cooperative-game} with the safety constraint \eqref{ineq:sufficient-condition} and its LP formulation in \eqref{opt: lp}. 
Suppose Assumption \ref{assumption: unichain} holds and let $\{X_{1}^{*}(s, a)\}, \{X_{2}^{*}(s, p) \}$ denote the optimal solutions of the LP problem \eqref{opt: lp}.
Then the optimal control and communication policies are given by
\begin{equation}
\label{eq: optimal-policy}
\mu_{\infty}^{m}(a)=\frac{X_{1}^{*}(s, a)}{\sum_{a \in A(s)}X_{1}^{*}(s, a)}, \quad \mu_{\infty}^{p}(p')=\frac{X_{2}^{*}(s, p')}{\sum_{p' \in \Omega_{p}}X_{2}^{*}(s, p')}.
\end{equation}
\end{theorem}
\begin{IEEEproof}
	By Proposition \ref{proposition: stationary-policy}, we only need to search for the stationary policy space of the optimal solutions to Problem \ref{problem: constrained-cooperative-game}. 
	Let $\overline{X}_{1}(s, a) \triangleq \lim_{T \rightarrow \infty}\frac{1}{T} \sum_{k=0}^{T}\mathbb{P}\{s_k=s, a_k=a\}$ denotes the average probability distribution of the state-action pair $(s, a)$, and $\overline{X}_{2}(s, p) \triangleq  \lim_{T \rightarrow \infty}\frac{1}{T} \sum_{k=0}^{T}\mathbb{P}\{s_k=s, p_k=p\}$ denotes the average probability distribution of the state-power pair $(s, p)$. Then, the joint cost $J(s_0, \mu^{\infty})$ under stationary policy $\mu^{\infty}$ is
	\begin{align*}
	\lim_{T \rightarrow \infty}\frac{1}{T} &\sum_{k=0}^{T}\mathbb{E}_{s_0}^{\mu^{\infty}} \big \{ c_{M}(s_k, a_k)+ \lambda c_{p}(p_k) \big\} \\
	&=\sum_{\substack{s \in S, a \in A(s) \\ p \in \Omega_{p}}} \lim_{T \rightarrow \infty}\frac{1}{T} \sum_{k=0}^{T}\mathbb{P}\{s_k=s, a_k=a\}c_{M}(s, a) \\
	& \qquad \qquad \qquad+\lambda c_{p}(p)\mathbb{P}\{ p_k = p\} \\
	&=\sum_{\substack{s \in S, a \in A(s) \\ p \in \Omega_{p}}} \overline{X}_{1}(s, a)c_{M}(s, a)+\lambda c_{p}(p)\overline{X}_{2}(s, p).
	\end{align*}
By the stationarity property of the control and communication policies (cf. Proposition \ref{proposition: stationary-policy} and  Assumption \ref{assumption: unichain}), there exists for a given $\mu_{\infty}=(\mu^{m}_{\infty}, \mu^{p}_{\infty})$  a unique stationary probability distribution $\{\overline{X}_{1}(s, a)\}$ and $\{\overline{X}_{2}(s, p)\}$ which satisfy the first two  constraints in the LP formulation \eqref{opt: lp} (cf. the proof of \cite[Theorem 4.3]{altman1999constrained} for details). 
It remains to show that the third constraint in \eqref{opt: lp} is equivalent to \eqref{ineq:sufficient-condition}. 
To do this, first note that $\overline{X}_{2}(s, p)=\mathbb{P}\{ p \vert s \} \lim_{T \rightarrow \infty}\frac{1}{T} \sum_{k=0}^{T} \mathbb{P}\{ s_{k}=s\}$ holds due to the use of stationary policy.
Furthermore, under Assumption \ref{assumption: unichain}, we know there exists a unique stationary probability distribution $\pi(s) \triangleq \mathbb{P}\{s\}$ and that the average probability distribution $\lim_{T \rightarrow \infty}\frac{1}{T} \sum_{k=0}^{T} \mathbb{P}\{ s_{k}=s\}$ converges to $\pi(s)$, i.e., $\pi(s)=\lim_{T \rightarrow \infty}\frac{1}{T} \sum_{k=0}^{T} \mathbb{P}\{ s_{k}=s\}$. 
These then imply that $\overline{X}_{2}(s, p)=\mathbb{P}\{p \vert s \} \pi(s)=\mathbb{P}\{s, p\} \triangleq \overline{\pi}(\overline{s})$ with $\overline{s}=(s, p)$. 
It is thus clear that the third constraint in \eqref{opt: lp} matches the safety condition in \eqref{ineq:sufficient-condition}. 
Hence, by the results in \cite[Theorem 4.3]{altman1999constrained}, it may then be concluded that  the optimal solutions of the LP problem \eqref{opt: lp} define the optimal policies  in \eqref{eq: optimal-policy}. 
\end{IEEEproof}
\begin{remark}
The equivalence between the optimal solutions of the constrained two-player game \eqref{eq:game} and the LP \eqref{opt: lp} lies in the simple structure of the safety condition in \eqref{ineq:sufficient-condition}, which can be represented as a linear combination of the decision variables.
This means that the optimization problem \eqref{eq:game}  will be much more challenging if the assigned safety constraints cannot be transformed into linear or convex constraints. 
Particularly as shown in the next section, if the stronger notion of ASAS needs to be enforced, the almost sure safety condition in \eqref{ineq: sufficient-safety} leads to a quadratic constraint. 
Section \ref{subsec: polynomial-opt} formulates a quadratic programming (QP) problem for solving the co-design problem to ensure the  ASAS property of NCS \eqref{eq: ncs}.  
\end{remark}

\subsection{Co-design of stochastic efficiency and ASAS: A quadratic programming approach}
\label{subsec: polynomial-opt}

Note that the sufficient conditions \eqref{ineq: sufficient-safety} for ASAS lead to polynomial constraints in the formulated optimization problem. 
In general, it is computationally hard to search for optimal policies over a history-dependent feasible space and polynomial constraints \eqref{ineq: sufficient-safety}. 
In what follows, we thus focus on finding optimal policies over the stationary policy space. 

Under the polynomial constraint \eqref{ineq: sufficient-safety}, the optimal solutions to Problem \ref{problem: constrained-cooperative-game} can be obtained by solving the following polynomial optimization problem for all $s \in S$ and $s' \in S$.
\begin{equation}
\begin{aligned}
& \underset{\{X(s, a, p)\}}{\min} 
& & \sum_{\substack{s \in S, a \in A(s) \\ p \in \Omega_{p}}} X(s, a, p)[c_{M}(s, a)+\lambda c_{p}(s, p)] \\
& \text{s. t.} &0=& \sum\limits_{a \in A(s')} X(s', a)-\sum\limits_{\substack{s \in S \\ a \in A(s)}} p(s' | s, a)X(s, a)\\
&&1=& \sum\limits_{\substack{s \in S, a \in A(s)\\ p \in \Omega_{p}}} X(s, a, p)\\
&&\overline{\eta}\triangleq &\frac{\eta-\lambda_{1} \rho}{\rho(\lambda_{0}-\lambda_{1})} \geq \sum\limits_{\substack{s' \in S \\ p' \in \Omega_{p}}} \frac{\theta(s', p')X(s, p')}{X(s)}\\
&&& \qquad\qquad\qquad\; \times\sum\limits_{a \in A(s)} \frac{p(s' | s, a)X(s, a)}{X(s)}.
\end{aligned}
\label{opt-polynomial}
\end{equation}
Note that the safety constraints in \eqref{opt-polynomial} may be transformed into equivalent quadratic constraints as follows
\begin{equation}
\overline{\eta}X^2(s)-\sum_{s' \in S}\bigg[\sum_{p' \in \Omega_{p}} \theta(s', p')X(s, p')\sum_{a \in A(s)}p(s' | s, a)X(s, a)\bigg] \geq 0
\label{ineq:quadratic}
\end{equation}
where $X(s) \neq 0, \forall s \in S$. 
Since $X(s)=\sum_{p' \in \Omega_{p}, a \in A(s)}X(s, a, p)$, $X(s, p')=\sum_{a \in A(s)}X(s, p', a)$ and $X(s, a)=\sum_{p' \in \Omega_{p}}=X(s, p', a)$, we introduce $\{x_{i}\}_{i=1}^{|S||A||\Omega_{p}|}$ as a set of variables which relabel the decision variables $\{X(s, a, p)\}_{s \in S, a \in A(s), p \in \Omega_{p}}$. 
For a given variable vector $x=[x_{1}, \ldots,x_{i}, \ldots, x_{|S||A||\Omega_{p}|}]$, it can be shown that there always exists a symmetric matrix $Q$ such that the constraints in \eqref{ineq:quadratic} can be rewritten in a compact form below
\begin{equation}
x^{T}Qx \geq 0.
\end{equation}
The elements of the symmetric matrix $Q$  are combinations of packet dropout probabilities $\{\theta(s, p): s \in S, p \in \Omega_{p}\}$ of the state-dependent wireless fading channel and transition probability $\{p(s' |s, a): s', s \in S, a \in A(s)\}$ of the MDP. 
It is clear that the optimization problem \eqref{opt-polynomial} is an LP problem with  quadratic constraints. 
In particular, the optimization problem \eqref{opt-polynomial} is convex if and only if the symmetric matrix $Q$ is negative semidefinite~(i.e., $Q \preceq 0$). 

\begin{theorem}
	Consider the co-design Problem \ref{problem: constrained-cooperative-game} with safety constraints \eqref{ineq: sufficient-safety}. 
Suppose $Q(\overline{\eta}, \theta, p)$  is a negative semi-definite matrix as defined by the quadratic inequalities in \eqref{ineq:quadratic}. 
Then the optimal stationary control and power policies of the co-design Problem \ref{problem: constrained-cooperative-game} can be obtained by solving the QP problem  in \eqref{opt-polynomial}. 
In particular, let $\{X^{*}(s, a, p)\}$ be the optimal solutions to the QP fomulation in \eqref{opt-polynomial}.
Then the optimal power and control policies are defined as follows.
	\begin{align}
	\mu^{m}_{s}(a)&=\frac{\sum_{p \in \Omega_{p}}X^{*}(s, a, p)}{\sum_{a \in A(s), p \in \Omega_{p}}X^{*}(s, a, p)},\\
	\mu^{p}_{s}(p)&=\frac{\sum_{a \in A(s)}X^{*}(s, a, p)}{\sum_{a \in A(s) p \in \Omega_{p}}X^{*}(s, a, p)}.
	\end{align}
\end{theorem}
\begin{IEEEproof}
The proof is similar to that of Theorem \ref{thm: lp} and thus omitted due to space limitation. 
\end{IEEEproof}
\begin{remark}
The convexity of the polynomial optimization problem in \eqref{opt-polynomial} solely depends on the negative semidefiniteness of the symmetric matrix $Q$.
As shown in inequality \eqref{ineq:quadratic}, $Q(\overline{\eta}, \theta, p)$ is a function of the conditions on the fading channel $\theta$ as well as the dynamics of both the MDP and the NCS as characterized by $p$ and $\overline{\eta}$, respectively. 
\end{remark}

The next section presents numerical examples  to illustrate the implementation of the proposed co-design framework in realistic NCS models.
\section{Simulation Example and Results}
\label{sec: simulation}
This section presents simulation results that were obtained when verifying the proposed co-design of safety~(Theorem \ref{thm-asas}) and optimal policies on dynamic model of a networked DC motor system and an autonomous vehicle~(e.g. forklift truck).
Comparison with more traditional design methods based on the separation principle (cf. \cite{gatsis2014optimal, rabi2016separated}) is also reported. 

\subsection{Stochastic stability of a networked DC motor}
In this simulation, a networked DC motor model is examined as an example of NCS due to its importance in manufacturing systems \cite{tipsuwan2003control}. 
Let $x=[\phi; \omega]$ be the state vector of the model consisting states $\phi$ and $\omega$ that describe the motor angular position and angular velocity, respectively. 
The dynamics of the DC motor satisfy the following linear differential equation model (cf. \cite{li2009state,li2009optimal}).
\begin{align}
\dot{x}(t)=\underbrace{\begin{bmatrix}
0 & 1 \\
1 & -217.4
\end{bmatrix}}_{A}x(t)+\underbrace{\begin{bmatrix}
0 \\
1669.5
\end{bmatrix}}_{B}u(t)
\end{align}
where $\dot{x}(t)=\frac{dx(t)}{dt}$ and $u(t)$ is the control input at time $t\in\mathbb{R}_{\geq 0}$.
Without loss of generality, let $x^{*}=0$ be the setpoint~(equilibrium) of the networked DC motor system\footnote{To be able to track non-zero reference signal, additional constraints must be imposed on the choice of those nonzero setpoints under the state feedback control law (cf. \cite{li2009eda} for more details).}. 
Let $T=0.3 s$ denotes the sampling time period for the system and the discretized DC motor model under such a sampling is 
\begin{align*}
x(k+1)=\underbrace{\begin{bmatrix}
1.0014 & 0.0046 \\
0.0046 & 0 
\end{bmatrix}}_{F}x(k)+\underbrace{\begin{bmatrix}
2.27 \\
7.6897
\end{bmatrix}}_{G} u(k)
\end{align*}
where $F=\exp(AT)$ and $G=\int_{0}^{T}\exp(A\tau)d\tau\,B$. 
At each discrete time instant $kT, k=0, 1, \ldots$, the communication system encodes the sampled state $x(k)$ into a single packet and transmit it through a wireless fading channel. 
The output of the fading channel is thus characterized by a random process $\{\gamma_k\}$ that is defined in \eqref{eq: SDDC} in which $\gamma_{k}=0$ represents the packet dropout event and $\gamma_{k}=1$ being the successful packet reception.  
Let $K=[-0.4055  -0.0024]$ be the chosen control gain vector such that $F+GK$ is Schur stable. 
The remote control signal\footnote{Rather than the zero-input strategy used in this paper, one can also adopt the hold-input control policy where the previous control input is taken if packet is lost, i.e., $ \hat{x}(k)=\gamma_k x(k)+(1-\gamma_k)\hat{x}(k-1), u(k)=K\hat{x}(k), \quad \hat{x}(0)=0$. 
However, as shown in \cite{schenato2009zero}, none of the two can be claimed superior to the other.} $u(k)$ is then constructed as $u(k)=\gamma_{k}Kx(k)$. 
Then the closed-loop of the networked DC motor system is a stochastic switched linear system $x(k+1)=(F+\gamma_{k}GK)x(k)$. 
Consider a Lyapunov function $V_{0}(x)=V_{1}(x)=x^{T}Px$ for this system with 
\begin{align*}
P=\begin{bmatrix}
6.5982  &   0.1143
\\
0.1143 &   0.0582
\end{bmatrix},
\end{align*}
such that we have $\lambda_{1}=0.1,\; \lambda_{0}=1.03$ and $\rho=1$.

\begin{figure}
	\centering
	\includegraphics[width=0.53\textwidth]{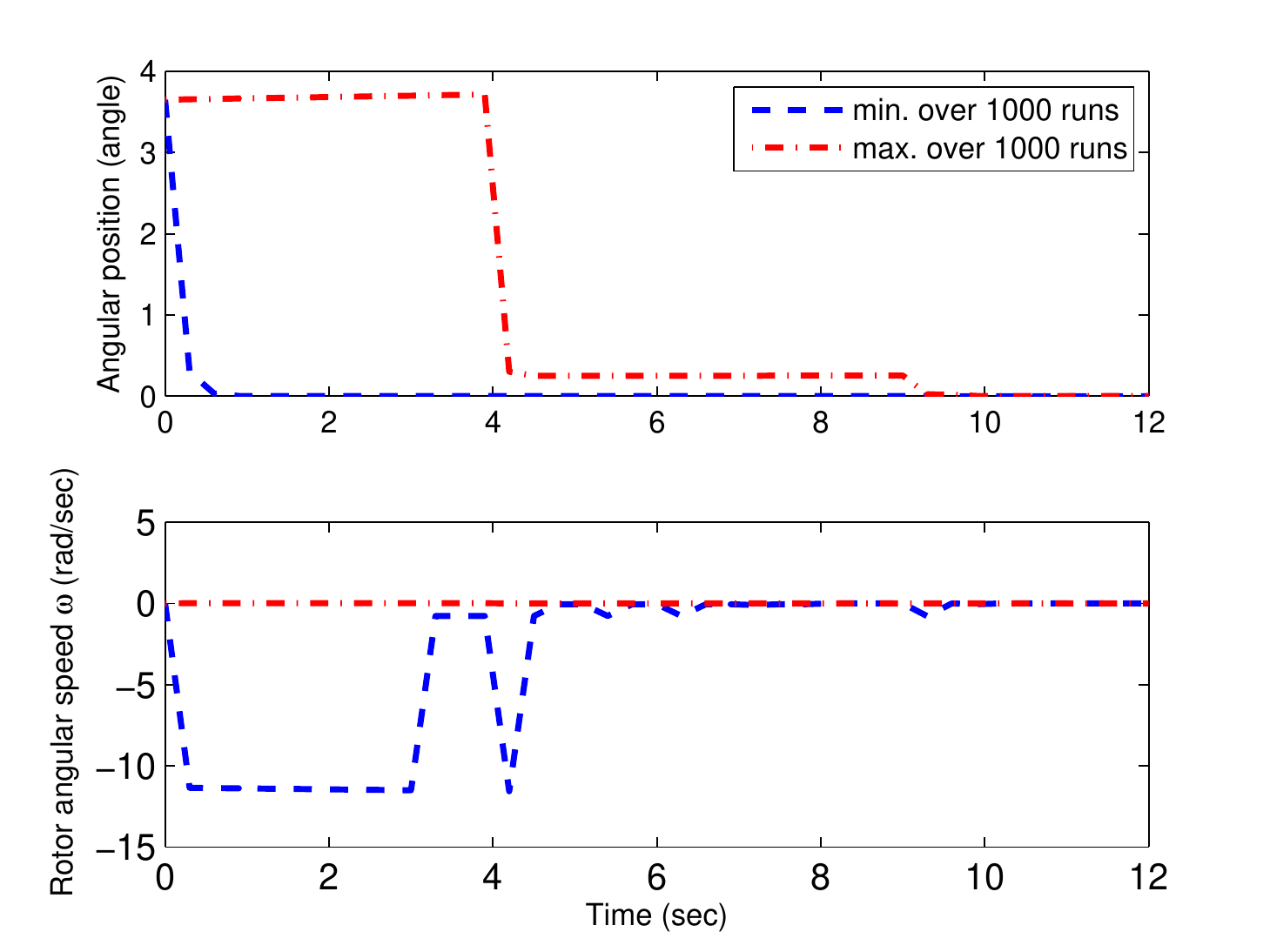}
	\caption{Maximum and minimum state values over $100$ samples.}
	\label{fig:as-as}
\end{figure}

With the selected system parameters $\lambda_{i}~(i=0,1)$ and $\rho$, the first part of the simulation is to verify the ASAS property under the condition in \eqref{ineq: sufficient-safety}.  The stochastic stability is evaluated using 100 runs of Monte Carlo simulation for $12$ seconds and with similar system parameters. 
Fig. \ref{fig:as-as} shows that the maximum~(marked by red dashed-dot line) and minimum~(marked by blue dashed line) value of the state trajectories over $100$ samples asymptotically converge to the setpoint $x^{*}=0$ as time increases. 
This proves that the networked DC motor system indeed posseses the ASAS property under the sufficient condition in \eqref{ineq: sufficient-safety}. 

\begin{figure}
	\centering
	\includegraphics[width=0.53\textwidth]{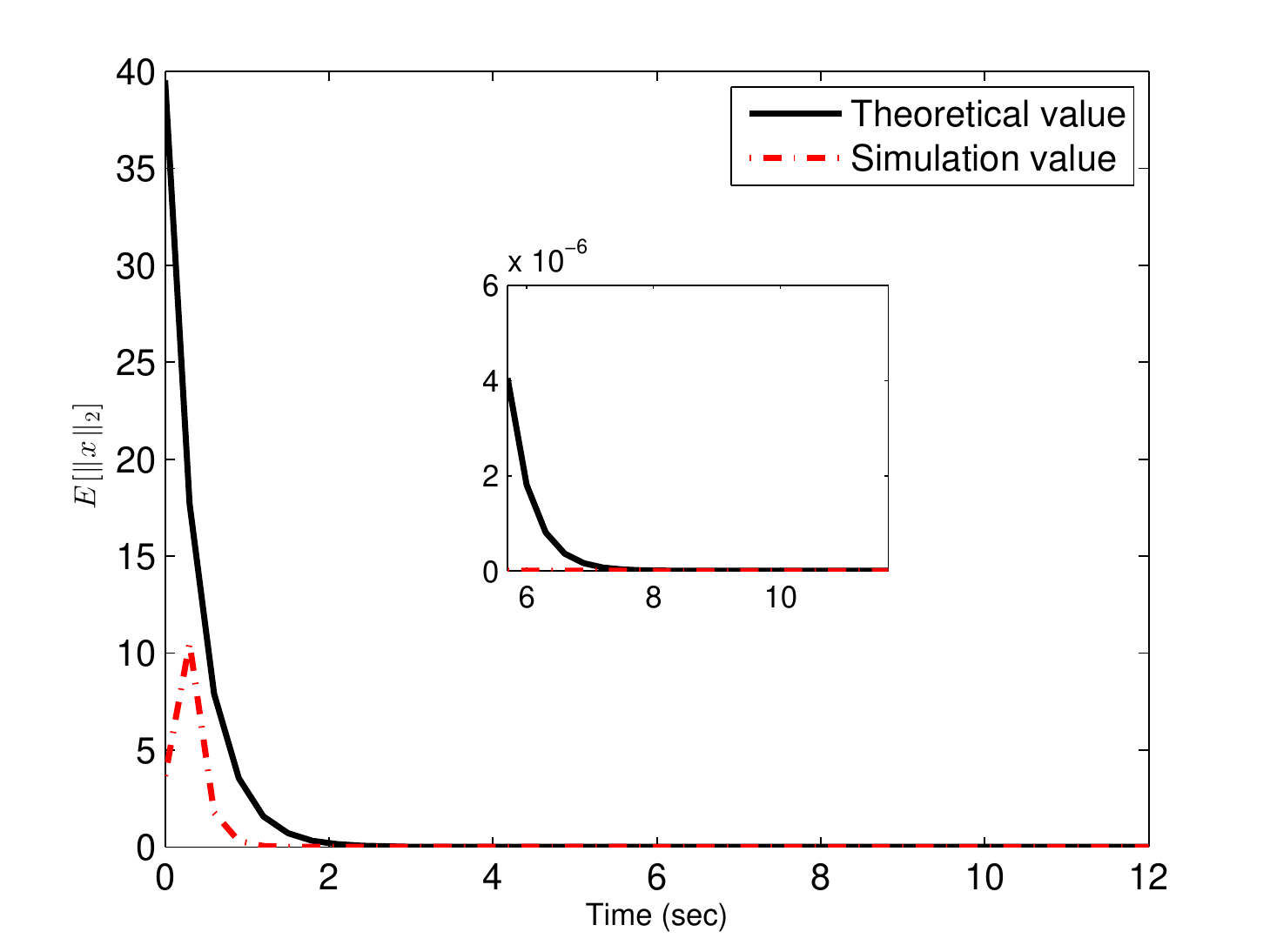}
	\caption{Comparison of theoretical and simulated expected state trajectories under a convergent rate of $\eta=0.1$.}
	\label{fig:as-es}
\end{figure}

Fig. \ref{fig:as-es} is used to verify the results of ASE presented in Theorem \ref{thm: safety-expectation}. As shown in Theorem \ref{thm: safety-expectation}, if the condition in \eqref{ineq: sufficient-safety} is satisfied with a convergence rate $\eta \in (0, 1)$  that is strictly less than one, a stronger notion of exponential safety in expectation can be ensured. 
To verify this result, let us choose $\eta=0.1$.
With this choice of $\eta$, Fig. \ref{fig:as-es} shows the comparison between the predicted theoretical value and the averaged value of $\|x\|_{2}$ for the 100 sample paths generated by Monte Carlo simulations. It can be seen on this figure that the result based on the Monte Carlo simulations is upper bounded by an exponential function that is predicted by the proposed theoretical results.

\subsection{MDP for an autonomous vehicle system}
As discussed previously, an MDP can be used to model the high-level transitions of autonomous vehicles (e.g. forklift truck, UGV) in manufacturing environments to accomplish some pick-and-deliver tasks. 
This simulation considers a three states MDP model to represents the three partitioned regions of an industrial workspace. 
At each state $s_i$, the autonomous vehicle can transit to one of the other two states by selecting appropriate actions $a_i$ with a predefined probability.  
The value of transition probabilities in MDP are summarized in Table \ref{table:tp-mdp}. 
For each state-action pair, a cost $c_M(\cdot)$ is defined to evaluate the system performance regarding the task accomplishment. 
These costs are also provided in Table \ref{table:tp-mdp}.

\begin{table}[!t]
	\caption{MDP Transition Probability~($P$) and Costs~($c_{M}$)}
	\label{table:tp-mdp}
	\centering
	\begin{tabular}{c|c c c||c}
		 & $s_1$ & $s_2$ & $s_3$ & $c_{M}(s, a)$ \\
		\hline
		$s_1, a_1$ & $1$ & $0$ & $0$ & $1$\\
		$s_1, b_1$ & $0.2$ & $0$ & $0.8$ & $2$\\
		\hline
		$s_2, a_2$ & $0.9$ & $0.1$ & $0.0$ & $2$\\
		$s_2, b_2$ & $0$ & $0.2$ & $0.8$ & $4$\\
		\hline
		$s_3, a_3$ & $0.1$ & $0.0$ & $0.9$ & $4$\\
		$s_3, b_3$ & $0$ & $0.8$ & $0.2$ & $2$\\
		\hline
	\end{tabular}
\end{table}
As stated in Section \ref{subsec:NNCS}, the SDDC model in \eqref{eq: SDDC} is used to model the impact that the MDP states have on the channel conditions~(packet dropout probability). For this example, the dropout probabilities $\theta(\cdot, \cdot)$ defined in \eqref{eq: SDDC} under the three MDP states~($s_1, s_2, s_3$) and the two power levels~($L, H$) are summarized in Table \ref{table:pr-sddc}. 
The power costs $c_{p}(p)$ are also specifed in Table \ref{table:pr-sddc}. 
The values of the dropout probability $\theta(\cdot, \cdot)$ are chosen in such a way that the simulation closely describes a situation where the levels of shadow fading decrease from the region labeled by $s_1$ to region of $s_3$. 
As indicated in Table \ref{table:pr-sddc}, the system can counteract the shadow fading phenomenon by selecting higher transmission power values.

\begin{figure}
		\centering
	\includegraphics[width=0.53\textwidth]{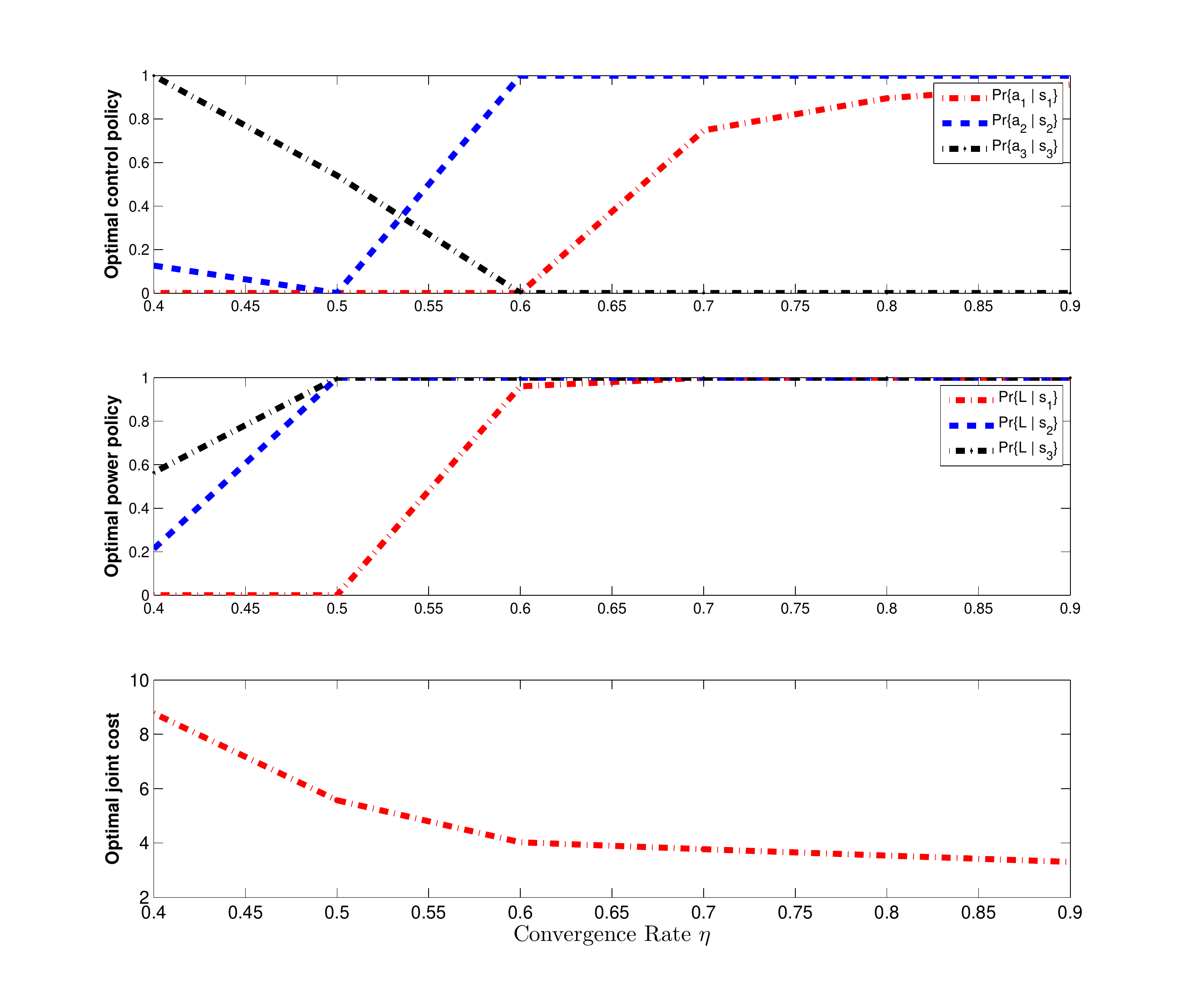}
	\caption{Optimal policies and costs under co-design framework over a wide range of convergence specifications for the networked DC motor system.}
	\label{fig:optimal-codesign}
\end{figure}

Fig. \ref{fig:optimal-codesign} shows the resulting optimal joint costs under the optimal co-design framework over a wide range of convergence rate (i.e. $0.4 \leq \eta \leq 0.9$). 
In this figure, the trade-off between safety and efficiency is demonstrated by adjustments that are made by optimal control and transmission power policies with regard to variation of convergence rates. 
In particular, when the convergence rate is small (which indicates high performance requirement for the networked DC motor system), the optimal control policy drives the forklift truck  away from the bad channel region $s_3$ which has a high shadow fading level.
At the same time, the optimal communication policy prefers to use a high transmission power level to ensure the attainment of the desired convergence rate.
In contrast, when the convergence rate increases, the optimal control and communication policies change to other directions to ensure efficient use of control and communication resources/energies, respectively. 
These simulation results thus clearly demonstrate the effectiveness of the proposed co-design framework to achieve both the safety and efficiency of industrial NCS applications.  

\begin{figure}
	\centering
	\includegraphics[width=0.53\textwidth]{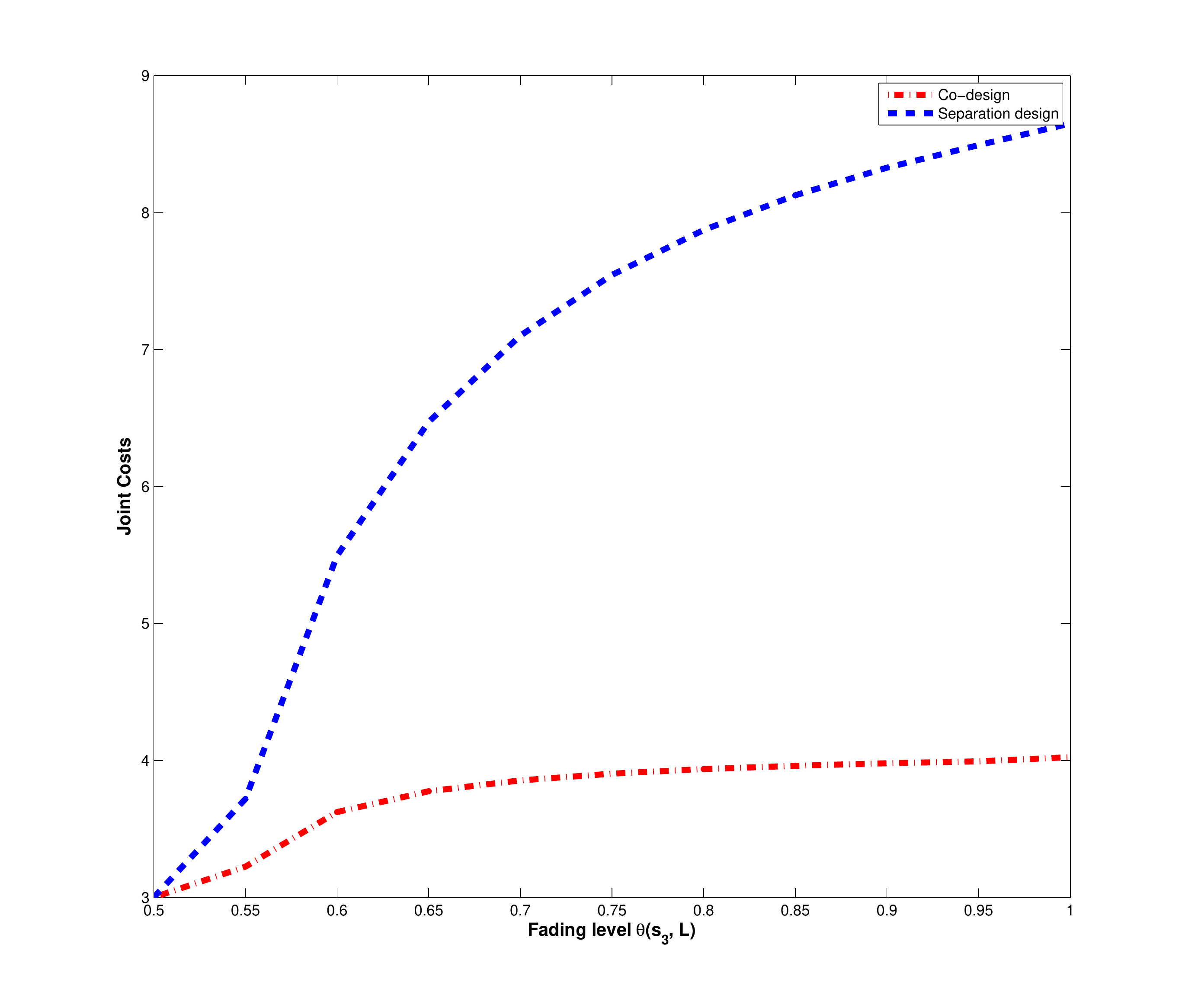}
	\caption{Performance comparison of co-design and separation design under different fading levels.}
	\label{fig:performance-comparison}
\end{figure}

Table \ref{table:comparison} and Fig. \ref{fig:performance-comparison} compares the results obtained by the proposed co-design strategy and that by the conventional separation method.
Specifically, this table compares the system performances that are achieved by both strategies under different convergence rates $\eta$ and fading levels $\theta$. 
As discussed in \cite{gatsis2014optimal, rabi2016separated}, the separation design approach generates the optimal policies under the assumption that the channel state~(packet dropout) is independent of the physical states. 
Fig. \ref{fig:performance-comparison} shows the comparison of the joint costs under the optimal co-design approach (plotted as a red dash-dot line) and the optimal separation design method (plotted as a blue dashed line) over a variety of fading levels within the range $0.5 \leq \theta(s_3, L) \leq 1$. 
It is clear from these plots that, while the optimal cost obtained by the co-design policy barely change as the fading level increases, the optimal cost obtained under the separation design policy  increases as the fading level increases. 
These results imply that the co-design policy leads to optimal performances for both communication and control systems that are more robust against shadow fading than those achieved by the separation policy. 
Table \ref{table:comparison} further provides numerical results which compare the joint costs achieved by both strategies under different $\theta$ and $\eta$. 
As shown by this table, the proposed co-design strategy outperforms the separation design in each scenario, either in the case of high fading levels with $\theta=0.85$ and $\theta=0.95$ or low fading levels with $\theta=0.65$ and $\theta=0.75$. 
It is worth noting that the separation design cannot find optimal policies that achieve the specified convergence rate $\eta$ for similar case of high shadow fading levels with $\theta=0.85$ and $\theta=0.95$. 
These results clearly demonstrate the necessity as well as the benefits of the proposed co-design strategy in achieving both the safety and efficiency requirements over the separation design method.  

\begin{table}[!t]
	\caption{Dropout Probability $\theta$ and Power Costs~($c_{p}$)}
	\label{table:pr-sddc}
	\centering
		\begin{tabular}{c|c c c||c}
			 & $s_1$ & $s_2$ & $s_3$ & $c_{p}(p)$ \\
			 \hline
			 L & $.9$ & $.5$ & $.4$ & $1$\\
			 H & $.4$ & $.3$ & $.2$ & $4$\\
			 \hline
		\end{tabular}
\end{table}

\begin{table}[!t]
\caption{Performance Comparison~($c_{M}+c_{p}$) under Co-design and Separation-design Framework}
\centering
\begin{tabular}{ll|l|l|l|l|}
\cline{3-6}
& & \multicolumn{4}{c|}{Convergence rate~$\eta$} \\
\hline
$\theta(s_1, L)$ & Policy & $0.4$ & $0.5$ & $0.6$ & $0.7$ \\ \hline
\multirow{2}{*}{$0.95$} & Co-design & $5.66$ & $5.66$ & $5.66$ & $5.66$ \\ 
& Sep-design & N/A & N/A & N/A & N/A \\ \hline
\multirow{2}{*}{$0.85$} & Co-design & $4.05$ & $4.05$ &$4.05$ & $4.05$ \\
& Sep-design & N/A & N/A & N/A & $7.33$ \\ \hline
\multirow{2}{*}{$0.75$} & Co-design & $2.99$ &$2.96$ & $2.90$ & $2.77$    \\
& Sep-design & $5.24$ & $5.04$ & $4.72$ & $4.12$ \\ \hline
\multirow{2}{*}{$0.65$} & Co-design & $2.73$ &$2.63$ & $2.45$ & $2.03$   \\
& Sep-design & $4.40$ & $4.0$ & $3.35$ & $2.09$ \\ \hline
\end{tabular}
\label{table:comparison}
\end{table}

\section{Conclusion}
\label{sec: conclusion}
This paper has examined the optimal co-design problem of industrial NCS under state-dependent correlated fading channels. 
The state-dependent property of the used wireless communication fading channels often arises from the movement of sizable machineries within the industrial workspace.
Such a fading has particularly been considered as a critical factor that could lead to safety-related issues in industrial automation systems. 
This paper explicitly characterized the property of such a dependency through the development of a state-dependent dropout channel model. 
Based on the developed model, this paper derived sufficient conditions which ensure the stochastic safety of industrial NCS in factory automation. 
Using the derived safety conditions, efficient co-design strategies which guarantee the attainment of both optimal control and communication performances are presented. 
Numerical simulation results on a model of heterogeneous industrial system which consists of a networked DC motor and an autonomous forklift truck were reported to verify the benefits and advantages of the proposed co-design framework over methods based on the conventional separation principle.

\bibliographystyle{IEEEtran}     
\bibliography{bibfile/TCNS-2018rev} 

\begin{thebibliography}{10}
\providecommand{\url}[1]{#1}
\csname url@samestyle\endcsname
\providecommand{\newblock}{\relax}
\providecommand{\bibinfo}[2]{#2}
\providecommand{\BIBentrySTDinterwordspacing}{\spaceskip=0pt\relax}
\providecommand{\BIBentryALTinterwordstretchfactor}{4}
\providecommand{\BIBentryALTinterwordspacing}{\spaceskip=\fontdimen2\font plus
\BIBentryALTinterwordstretchfactor\fontdimen3\font minus
  \fontdimen4\font\relax}
\providecommand{\BIBforeignlanguage}[2]{{%
\expandafter\ifx\csname l@#1\endcsname\relax
\typeout{** WARNING: IEEEtran.bst: No hyphenation pattern has been}%
\typeout{** loaded for the language `#1'. Using the pattern for}%
\typeout{** the default language instead.}%
\else
\language=\csname l@#1\endcsname
\fi
#2}}
\providecommand{\BIBdecl}{\relax}
\BIBdecl

\bibitem{song2008wirelesshart}
J.~Song \emph{et~al.}, ``Wirelesshart: Applying wireless technology in
  real-time industrial process control,'' in \emph{Proc. IEEE RTAS}, St. Louis,
  MO, USA, 2008, pp. 377--386.

\bibitem{gungor2009industrial}
V.~C. Gungor and G.~P. Hancke, ``Industrial wireless sensor networks:
  Challenges, design principles, and technical approaches,'' \emph{IEEE Trans.
  Ind. Electron.}, vol.~56, no.~10, pp. 4258--4265, 2009.

\bibitem{agrawal2014long}
P.~Agrawal, A.~Ahl{\'e}n, T.~Olofsson, and M.~Gidlund, ``Long term channel
  characterization for energy efficient transmission in industrial
  environments,'' \emph{IEEE Trans. Commun.}, vol.~62, no.~8, pp. 3004--3014,
  2014.

\bibitem{aakerberg2011future}
J.~{\AA}kerberg, M.~Gidlund, and M.~Bj{\"o}rkman, ``Future research challenges
  in wireless sensor and actuator networks targeting industrial automation,''
  in \emph{Proc. IEEE INDIN}, Lisbon, Portugal, 2011, pp. 410--415.

\bibitem{willig2005wireless}
A.~Willig, K.~Matheus, and A.~Wolisz, ``Wireless technology in industrial
  networks,'' \emph{Proc. IEEE}, vol.~93, no.~6, pp. 1130--1151, 2005.

\bibitem{willig2008recent}
A.~Willig, ``Recent and emerging topics in wireless industrial communications:
  A selection,'' \emph{IEEE Trans. Ind. Informat.}, vol.~4, no.~2, pp.
  102--124, 2008.

\bibitem{kashiwagi2010time}
I.~Kashiwagi, T.~Taga, and T.~Imai, ``Time-varying path-shadowing model for
  indoor populated environments,'' \emph{IEEE Trans. Veh. Technol.}, vol.~59,
  no.~1, pp. 16--28, 2010.

\bibitem{quevedo2013state}
D.~E. Quevedo, A.~Ahlen, and K.~H. Johansson, ``State estimation over sensor
  networks with correlated wireless fading channels,'' \emph{IEEE Trans. Autom.
  Control}, vol.~58, no.~3, pp. 581--593, 2013.

\bibitem{olofsson2016modeling}
T.~Olofsson, A.~Ahl{\'e}n, and M.~Gidlund, ``Modeling of the fading statistics
  of wireless sensor network channels in industrial environments,'' \emph{IEEE
  Trans. Signal Process.}, vol.~64, no.~12, pp. 3021--3034, 2016.

\bibitem{tse2005fundamentals}
D.~Tse and P.~Viswanath, \emph{Fundamentals of Wireless Communication}.\hskip
  1em plus 0.5em minus 0.4em\relax Cambridge University Press, 2005.

\bibitem{hu2017co}
B.~Hu, Y.~Wang, P.~Orlik, T.~Koike-Akino, and J.~Guo, ``Co-design of safe and
  efficient networked control systems in factory automation with
  state-dependent wireless fading channels,'' \emph{arXiv:1708.06468}, 2017.

\bibitem{gatsis2014optimal}
K.~Gatsis, A.~Ribeiro, and G.~J. Pappas, ``Optimal power management in wireless
  control systems,'' \emph{IEEE Trans. Autom. Control}, vol.~59, no.~6, pp.
  1495--1510, 2014.

\bibitem{ren2018infinite}
X.~Ren, J.~Wu, K.~H. Johansson, G.~Shi, and L.~Shi, ``Infinite horizon optimal
  transmission power control for remote state estimation over fading
  channels,'' \emph{IEEE Trans. Autom. Control}, vol.~63, no.~1, pp. 85--100,
  2018.

\bibitem{wen2018transmission}
S.~Wen, G.~Guo, B.~Chen, and X.~Gao, ``Transmission power scheduling and
  control co-design for wireless sensor networks,'' \emph{Inform. Sciences},
  vol. 442, pp. 114--127, 2018.

\bibitem{rabi2016separated}
M.~Rabi, C.~Ramesh, and K.~H. Johansson, ``Separated design of encoder and
  controller for networked linear quadratic optimal control,'' \emph{SIAM J.
  Control Optim.}, vol.~54, no.~2, pp. 662--689, 2016.

\bibitem{leong2017event}
A.~S. Leong, D.~E. Quevedo, T.~Tanaka, S.~Dey, and A.~Ahl{\'e}n, ``Event-based
  transmission scheduling and {LQG} control over a packet dropping link,''
  \emph{IFAC-PapersOnLine}, vol.~50, no.~1, pp. 8945--8950, 2017.

\bibitem{peng2013event}
C.~Peng and T.~C. Yang, ``Event-triggered communication and ${H}_\infty$
  control co-design for networked control systems,'' \emph{Automatica},
  vol.~49, no.~5, pp. 1326--1332, 2013.

\bibitem{peng2013novel}
C.~Peng and Q.-L. Han, ``A novel event-triggered transmission scheme and
  $\mathcal{L}_2$ control co-design for sampled-data control systems,''
  \emph{IEEE Trans. Autom. Control}, vol.~58, no.~10, pp. 2620--2626, 2013.

\bibitem{zhao2008integrated}
Y.~Zhao, G.~Liu, and D.~Rees, ``Integrated predictive control and scheduling
  co-design for networked control systems,'' \emph{IET Control Theory A.},
  vol.~2, no.~1, pp. 7--15, 2008.

\bibitem{varma2016energy}
V.~S. Varma and R.~Postoyan, ``Energy efficient time-triggered control over
  wireless sensor/actuator networks,'' in \emph{Proc. IEEE CDC}, LasVegas, NV,
  USA, 2016, pp. 2727--2732.

\bibitem{lyu2017co}
L.~Lyu, C.~Chen, C.~Hua, S.~Zhu, and X.~Guan, ``Co-design of stabilisation and
  transmission scheduling for wireless control systems,'' \emph{IET Control
  Theory A.}, vol.~11, no.~11, pp. 1767--1778, 2017.

\bibitem{peters2016controller}
E.~G. Peters, D.~E. Quevedo, and M.~Fu, ``Controller and scheduler codesign for
  feedback control over {IEEE} 802.15. 4 networks,'' \emph{IEEE Trans. Control
  Syst. Technol.}, vol.~24, no.~6, 2016.

\bibitem{eriksson2016long}
M.~Eriksson and T.~Olofsson, ``On long-term statistical dependences in channel
  gains for fixed wireless links in factories,'' \emph{IEEE Trans. Commun.},
  vol.~64, no.~7, pp. 3078--3091, 2016.

\bibitem{quevedo2014power}
D.~E. Quevedo, J.~{\O}stergaard, and A.~Ahlen, ``Power control and coding
  formulation for state estimation with wireless sensors,'' \emph{IEEE Trans.
  Control Syst. Technol.}, vol.~22, no.~2, pp. 413--427, 2014.

\bibitem{leong2016network}
A.~S. Leong, D.~E. Quevedo, A.~Ahl{\'e}n, and K.~H. Johansson, ``On network
  topology reconfiguration for remote state estimation,'' \emph{IEEE Trans.
  Autom. Control}, vol.~61, no.~12, pp. 3842--3856, 2016.

\bibitem{de2006use}
F.~De~Pellegrini, D.~Miorandi, S.~Vitturi, and A.~Zanella, ``On the use of
  wireless networks at low level of factory automation systems,'' \emph{IEEE
  Trans. Ind. Informat.}, vol.~2, no.~2, pp. 129--143, 2006.

\bibitem{kozin1969survey}
F.~Kozin, ``A survey of stability of stochastic systems,'' \emph{Automatica},
  vol.~5, no.~1, pp. 95--112, 1969.

\bibitem{puterman1994markov}
M.~L. Puterman, \emph{Markov Decision Processes: Discrete Stochastic Dynamic
  Programming}.\hskip 1em plus 0.5em minus 0.4em\relax John Wiley \& Sons,
  1994, vol.~10.

\bibitem{anthonisse1977exponential}
J.~M. Anthonisse and H.~Tijms, ``Exponential convergence of products of
  stochastic matrices,'' \emph{J. Math. Anal. Appl.}, vol.~59, no.~2, pp.
  360--364, 1977.

\bibitem{jiang2001input}
Z.-P. Jiang and Y.~Wang, ``Input-to-state stability for discrete-time nonlinear
  systems,'' \emph{Automatica}, vol.~37, no.~6, pp. 857--869, 2001.

\bibitem{liberzon2012switching}
D.~Liberzon, \emph{Switching in Systems and Control}.\hskip 1em plus 0.5em
  minus 0.4em\relax Birkh\"{a}user, 2012.

\bibitem{chatterjee2011stabilizing}
D.~Chatterjee and D.~Liberzon, ``Stabilizing randomly switched systems,''
  \emph{SIAM J. Control Optim.}, vol.~49, no.~5, pp. 2008--2031, 2011.

\bibitem{altman1999constrained}
E.~Altman, \emph{Constrained Markov Decision Processes}.\hskip 1em plus 0.5em
  minus 0.4em\relax CRC Press, 1999.

\bibitem{tipsuwan2003control}
Y.~Tipsuwan and M.-Y. Chow, ``Control methodologies in networked control
  systems,'' \emph{Control Eng. Pract.}, vol.~11, no.~10, pp. 1099--1111, 2003.

\bibitem{li2009state}
H.~Li, M.-Y. Chow, and Z.~Sun, ``State feedback stabilisation of networked
  control systems,'' \emph{IET Control Theory A.}, vol.~3, no.~7, pp. 929--940,
  2009.

\bibitem{li2009optimal}
------, ``Optimal stabilizing gain selection for networked control systems with
  time delays and packet losses,'' \emph{IEEE Trans. Control Syst. Technol.},
  vol.~17, no.~5, pp. 1154--1162, 2009.

\bibitem{li2009eda}
------, ``{EDA}-based speed control of a networked {DC} motor system with time
  delays and packet losses,'' \emph{IEEE Trans. Ind. Electron.}, vol.~56,
  no.~5, pp. 1727--1735, 2009.

\bibitem{schenato2009zero}
L.~Schenato, ``To zero or to hold control inputs with lossy links?'' \emph{IEEE
  Trans. Autom. Control}, vol.~54, no.~5, pp. 1093--1099, 2009.

\bibitem{ozekici1997markov}
S.~{\"O}zekici, ``Markov modulated {B}ernoulli process,'' \emph{Math. Method
  Oper. Res.}, vol.~45, no.~3, pp. 311--324, 1997.

\end{thebibliography}
\appendices
\section{Proof}
\label{appendix:proof}

\begin{IEEEproof}[Proof of Theorem \ref{thm: safety-expectation}]
	Consider the randomly switched nonlinear system in \eqref{eq: randomly-switched} with $w=0$ and the randomly switching signal $\{\gamma_{k}\}$ as defined by SDDC channel model in \eqref{eq: SDDC}.
Suppose that Assumption \ref{assumption: multiple-lyapunov} is satisfied by this system.
One then has
\begin{align*}
V_{\gamma_{k+1}}(x_{k+1}) & \overset{(a)}{\leq} \varrho V_{\gamma_k}(x_{k+1}) \overset{(b)}{\leq} \varrho \lambda_{\gamma_k}V_{\gamma_k}(x_k) \\
& \leq \varrho^{k+1} \prod_{\ell=0}^{k}\lambda_{\gamma_{\ell}}V_{\gamma_0}(x_0) \overset{(c)}{\leq} \varrho^{k+1}\prod_{\ell=0}^{k}\lambda_{\gamma_{\ell}}\alpha_{2}(|x_0|).
\end{align*}
The inequalities $(a), (b), (c)$ above hold by properties (V3), (V2) and (V1), respectively, in Assumption \ref{assumption: multiple-lyapunov}. 
Recall that the random switching signal $\{\gamma_{k}\}$ is a Markov-modulated Bernoulli process that is conditionally independent over time for the given joint-state $\overline{s}_k=(s_k, p_k) \in S \times \Omega_{p}$ \cite{ozekici1997markov}.
Let $\overline{s}_{0:k}=\overline{s}_{0} \ldots \overline{s}_{k}$ denotes the sample paths of the joint state up to time $k$ and $\overline{S}_{0:k}$ denotes the set of all the possible sample paths $\overline{s}_k$.
Correspondingly, let $\gamma_{0:k}$ denotes the sample paths of the randomly switching signal and $\Gamma_{0:k}$ denotes the set of all possible paths $\gamma_{0: k}$. 
By the Markovian property result obtained in Proposition \ref{proposition: transition-probability}, one knows that both of the sets $\overline{S}_{0:k}$ and $\Gamma_{0: k}$ are measurable.
As such, the expectation of the Lyapunov function $V_{\gamma_{k+1}}(x_{k+1})$ may be written as
\begin{align}
\mathbb{E}\big[V_{\gamma_{k+1}}&(x_{k+1})\big] \nonumber\\
 &\leq \varrho^{k+1} \mathbb{E}_{\overline{S}_{0:k}}\bigg[ \mathbb{E}_{\Gamma_{0: k}}\big[\prod_{\ell=0}^{k} \lambda_{\gamma_{\ell}} \big\vert \overline{s}_{0:k} \big]\bigg]\alpha_{2}(|x_0|) \nonumber\\
&\leq \varrho^{k+1}\mathbb{E}_{\overline{S}_{0:k}}\bigg[\prod_{\ell=0}^{k}\mathbb{E}[\lambda_{\gamma_{\ell}} \vert \overline{s}_{\ell}] \bigg]\alpha_{2}(|x_0|) \nonumber\\
&\leq 
\varrho^{k+1}\mathbb{E}_{\overline{S}_{0:k}}\bigg[\prod_{\ell=0}^{k}\big[(\lambda_{0}-\lambda_{1})\theta(\overline{s}_{\ell})+\lambda_{1}\big] \bigg]\alpha_{2}(|x_0|) \nonumber\\
&\leq \varrho^{k+1} \prod_{\ell=0}^{k} \bigg[ (\lambda_{0}-\lambda_{1})\theta^{T}\bm{\xi}_{\ell}+\lambda_{1}\bigg]\alpha_{2}(|x_0|)
\label{ineq: E-V}
\end{align}
where $\bm{\xi}_{\ell}=[{\rm Pr}\{\overline{s}_{\ell}=\overline{s}\}]_{|S||\Omega_{p}| \times 1}$ is the vector of probability distributions of the joint state $\overline{s}_{\ell}$ at time $\ell$ and $\theta=[\theta(\overline{s})]_{|S||\Omega_{p}| \times 1}$ is the corresponding vector of probabilities of the packet loss for all joint states $\overline{s} \in S \times \Omega_{p}$.
 For a selected control policy $\mu=\{\mu_{k}\}_{k=1}^{\infty}$, let $\overline{P}_{k}:=\overline{P}(\mu_{k})$ denotes the transition probability matrix of a Markov chain that is induced by the control policy $\mu$.
Then for an initial distribution $\bm{\xi}_{0}$, we have $\bm{\xi}_{k}=\prod_{i=1}^{k}\overline{P}_{i}\bm{\xi}_{0}$.
Thus, the inequality in  \eqref{ineq: E-V} may be rewritten as
\begin{align}
&\mathbb{E}\Big[V_{\gamma_{k+1}}(x_{k+1})\Big] \nonumber\\
&\leq \varrho^{k+1}\prod_{\ell=0}^{k}  \Big[(\lambda_{0}-\lambda_{1})\theta^{T}\prod_{i=1}^{k}\overline{P}_{i}\bm{\xi}_{0}+\lambda_{1}    \Big]\alpha_{2}(|x_0|) \nonumber \\
&\leq \prod_{\ell=0}^{k} \varrho \Big[(\lambda_{0}-\lambda_{1})\theta^{T}([\overline{\pi}_{i}+\alpha^{[\ell/v_{\ell}]}]_{ij})\bm{\xi}_{0}+\lambda_{1} \Big] \alpha_{2}(|x_0|)  \label{ineq: exp-convergece-P} \\
& = \prod_{\ell=0}^{k}\varrho\Big[(\lambda_{0}-\lambda_{1})\theta^{T}(\overline{\pi} \bm{1}^{T}+\bm{11}^{T}\alpha^{[\ell/v_{\ell}]}])\bm{\xi}_{0}+\lambda_{1} \Big] \alpha_{2}(|x_0|) \nonumber \\
&=\prod_{\ell=0}^{k}\varrho\Big[(\lambda_{0}-\lambda_{1})\theta^{T}(\overline{\pi} \underbrace{\bm{1}^{T}\bm{\xi}_{0}}_{=1}+\bm{11}^{T}\bm{\xi}_{0}\alpha^{[\ell/v_{\ell}]}])+\lambda_{1} \Big] \alpha_{2}(|x_0|) \nonumber \\
&=\prod_{\ell=0}^{k}\varrho\Big[(\lambda_{0}-\lambda_{1})\theta^{T}(\overline{\pi}+ \alpha^{\ell/v_{\ell}}\bm{1})+\lambda_{1}\Big]\alpha_{2}(|x_0|) \label{ineq:E-V-2}
\end{align}
where $\bm{1}=[1, 1, \ldots, 1]_{|S||\Omega_{p}| \times 1}$ is a column vector with all ements equal 1.
The inequality in \eqref{ineq: exp-convergece-P} holds due to the condition \eqref{ineq-convergence} in Theorem \ref{thm:old} and $\overline{\pi}(\mu)$ is the stationary distribution over the joint state $\overline{S}$ whose value depends on the selected control policy $\mu$. 

Now let $A_{\ell}:=\varrho(\lambda_{0}-\lambda_{1})\theta^{T}\bm{1}\alpha^{[\ell / v_{\ell}]}$ and $B:=\varrho(\lambda_{0}-\lambda_{1})\theta^{T}\overline{\pi}+ \varrho\lambda_{1}$.
Then, inequality \eqref{ineq:E-V-2} can be further rewritten as $\mathbb{E}\Big[V_{\gamma_{k+1}}(x_{k+1})\Big] \leq \prod_{\ell=0}^{k}(A_{\ell}+B)$ and its convergence  can be analyzed by examining dynamics of the series $\{A_{\ell}+B\}_{\ell=0}^{k}$ where $A_{\ell}$ and $B$ denote the transient and steady state parts, respectively. 
Since the infinite series of the transient part $\{A_{\ell}\}_{\ell=0}^{\infty}$ converges and $\lim_{\ell \rightarrow \infty} A_{\ell}=0$ due to the exponential convergence of $\alpha^{[\ell/v_{\ell}]}$ for $0 \leq \alpha < 1$, it is clear that the sufficient condition to ensure the whole infinite series produces $\prod_{\ell=0}^{\infty}(A_{\ell}+B)\leq 1$ is that $0 < B < 1$.
These lead to a condition which can be stated as follows.
\begin{align*}
0 \leq \varrho(\lambda_{0}-\lambda_{1})\theta^{T}\overline{\pi}+ \varrho\lambda_{1} < 1 \Leftrightarrow \theta^{T}\overline{\pi} < \frac{1-\varrho\lambda_{1}}{\varrho(\lambda_{0}-\lambda_{1})}.
\end{align*}
Since $\varrho\lambda_{1} < 1$ and $\lambda_{0} > \lambda_{1}$ hold true by property (V4) in Assumption \ref{assumption: multiple-lyapunov} and by the definitions of $\lambda_{0}$ and $\lambda_{1}$, the right hand side of inequality \eqref{ineq:sufficient-condition} is guaranteed to be positive as claimed in the theorem. 
The proof is thus completed. 
\end{IEEEproof}

\begin{IEEEproof}[Proof of Theorem \ref{thm-asas}]
Consider the randomly switched nonlinear system \eqref{eq: randomly-switched} without external disturbance ($w=0$) and a Lyapunov function $V_{\gamma_{k+1}}(x_{k+1})$.
Let $\mathbbm{1}_{A}: X \rightarrow \{0, 1\}$ denotes an indicator function of a subset $A$ of a set $X$, i.e., $\mathbbm{1}_{A}(x)=1$ if $x \in A$ and $\mathbbm{1}_{A}(x)=0$ otherwise.
We then have
\begin{align}
\mathbb{E}\Big[V&(x_{k+1}) \mathbbm{1}_{\{\overline{s}_{k+1}=\overline{s}\}} \Big] = \sum_{\ell=0}^{1} \mathbb{E} \Big[V_{\gamma_{k+1}}(x_{k+1}) \mathbbm{1}_{\{\overline{s}_{k+1}=\overline{s}\}} \mathbbm{1}_{\{\gamma_{k+1}=\ell\}} \Big]
\nonumber \\
&\leq \sum_{\ell=0}^{1}\mathbb{E}\Big[ \lambda_{\gamma_{k+1}} V_{\gamma_{k+1}}(x_k)  \mathbbm{1}_{\{\overline{s}_{k+1}=\overline{s}\}} \mathbbm{1}_{\{\gamma_{k+1}=\ell\}}) \Big] \nonumber \\
& \leq \sum_{\ell=0}^{1}\mathbb{E}\Big[ \lambda_{\gamma_{k+1}} \varrho V_{\gamma_{k}}(x_k)  \mathbbm{1}_{\{\overline{s}_{k+1}=\overline{s}\}} \mathbbm{1}_{\{\gamma_{k+1}=\ell\}}) \Big] \nonumber \\
&\leq \varrho\underbrace{\big(\lambda_{0}\theta(\overline{s})+(1-\theta(\overline{s}))\lambda_{1}\big)}_{\overline{\theta}(\overline{s})}\mathbb{E}\Big[ V_{\gamma_k}(x_k) \mathbbm{1}_{\{\overline{s}_{k+1}=\overline{s}\}}\Big] \nonumber \\
&= \varrho \overline{\theta}(\overline{s})\sum_{\overline{s}' \in \overline{S}}\mathbb{E}\Big[ V_{\gamma_k}(x_k) \mathbbm{1}_{\{\overline{s}_{k+1}=\overline{s}\}} \mathbbm{1}_{\{\overline{s}_{k=\overline{s}'}\}}\Big] \nonumber \\
&=\varrho \overline{\theta}(\overline{s})\sum_{\overline{s}' \in \overline{S}} \overline{P}(\overline{s}', \overline{s}) \mathbb{E} \Big[V_{\gamma_k}(x_k)\mathbbm{1}_{\{\overline{s}_{k}=\overline{s}' \}} \Big]
\label{ineq:E-V-s}
\end{align}
Let $\overline{V}_{k}=\big[\overline{V}_{k}(\overline{s})\big]_{\overline{s} \in \overline{S}}$ be a vector with $\overline{V}_{k}(\overline{s})\triangleq \mathbb{E}\big[V(x_{k})\mathbbm{1}_{\{\overline{s}_{k}=\overline{s}\}}\big]$.
Then the vectorial form of inequality \eqref{ineq:E-V-s} may be written as
\begin{align}
\overline{V}_{k+1} \leq \varrho \text{diag}\big(\overline{\theta}(\overline{s}_{1}), \ldots, \overline{\theta}(\overline{s}_{NM})\big) \overline{P} \overline{V}_{k}
\label{ineq:bar-V-sys}
\end{align} 
where $\overline{P}$ is the transition matrix of the joint state as defined in \eqref{eq: P_k}. 
Furthermore, since $\mathbb{E}\big[V(x_{k+1})\big]=\|\overline{V}_{k+1}\|_{1}$, we have
\begin{align*}
\mathbb{E}\big[V(x_{k+1})\big] \leq \varrho \|\text{diag}\big(\overline{\theta}(\overline{s}_{1}), \ldots, \overline{\theta}(\overline{s}_{NM})\big) \overline{P}\|_{1} \mathbb{E}\big[V(x_{k})\big]
\end{align*} 
It is clear that, if $\varrho \|\text{diag}\big(\overline{\theta}(\overline{s}_{1}), \ldots, \overline{\theta}(\overline{s}_{NM})\big) \overline{P}\|_{1} < 1$ (which is equivalent to condition \eqref{ineq: sufficient-safety}) holds, there exists a real $\eta$ with $0 \leq \eta < 1$ such that $\mathbb{E}\big[V(x_{k+1})\big] \leq \eta^{k+1} \alpha_{2}(|x_0|)$.

Next, we prove the ASAS property.
First, for a selected time $k \in \mathbb{Z}_{\geq 0}$, let $\underline{k} \leq k \leq \overline{k}$, $\forall r \geq 0$.
Consider the probability below
\begin{align*}
\mathbb{P}&\big\{ \sup_{\underline{k} \leq k \leq \overline{k}} V(x_k) \geq r \big\} \leq \frac{\mathbb{E}\big[ \sup_{\underline{k} \leq k \leq \overline{k}} V(x_k) \big]}{r} \\
&\qquad\qquad\qquad\leq \frac{\mathbb{E}\big[ \sum_{k=\underline{k}}^{\overline{k}} V(x_k) \big]}{r}
\leq \frac{\sum_{k=\underline{k}}^{\overline{k}} \mathbb{E}\big[ V(x_k) \big]}{r} \\
&\qquad\qquad\qquad\leq \sum_{k=\underline{k}}^{\overline{k}} \eta^{k} \frac{\alpha_{2}(|x_0|)}{r} \leq \eta^{\underline{k}}\frac{1-\eta^{\overline{k}-\underline{k}}}{1-\eta}\frac{\alpha_{2}(|x_0|)}{r}
\end{align*}
Letting $\overline{k} \rightarrow \infty$, one has that 
\begin{align}
\label{ineq: sup-V}
\mathbb{P}\big\{ \sup_{\underline{k} \leq k} V(x_k) \geq r \big\} \leq \frac{\eta^{\underline{k}}}{1-\eta}\frac{\alpha_{2}(|x_0|)}{r}.
\end{align}
Now since $\alpha_{1}(|x|) \geq r \Rightarrow V(x) \geq r$ holds due to property (V1) in Assumption \ref{assumption: multiple-lyapunov}, one has that $\mathbb{P}\big\{ \sup_{\underline{k} \leq k} |x| \geq \alpha_{1}^{-1}(r) \big\} \leq \mathbb{P}\big\{ \sup_{\underline{k} \leq k} V(x_k) \geq r \big\}$. 
Let $r' \triangleq \alpha_{1}^{-1}(r)$, one has that
\begin{align*}
\mathbb{P}\big\{ \sup_{\underline{k} \leq k} |x_k| \geq r \big\} \leq \frac{\eta^{\underline{k}}}{1-\eta}\frac{\alpha_{2}(|x_0|)}{r}\triangleq \xi(|x_0|, \underline{k}, r).
\end{align*}
Summing the above probability from $\underline{k}=1$ to $\underline{k}=\infty$ gives
\begin{align}
\sum_{\underline{k}=1}^{\infty}\mathbb{P}\big\{ \sup_{\underline{k} \leq k} |x_k| \geq r \big\} \leq  \sum_{\underline{k}=1}^{\infty} \frac{\eta^{\underline{k}}}{1-\eta}\frac{\alpha_{2}(|x_0|)}{r} < \infty \label{ineq: borel-cantelli-lemma}
\end{align}
Since $\eta \in [0, 1)$, the above infinite sum is finite.
By Borel-Cantelli Lemma, \eqref{ineq: borel-cantelli-lemma} becomes
$
\mathbb{P}\big\{ \lim_{\underline{k} \rightarrow \infty}\sup_{\underline{k} \leq k} |x_k| \geq r \big\}=0.
$
The proof is complete.  
\end{IEEEproof}

\begin{IEEEproof}[Proof of Theorem \ref{thm: practical-safety}]
Using the same arguments (under Assumption \ref{assumption: multiple-lyapunov}) as in the proof of Theorem \ref{thm: practical-safety}, the Lyapunov function $V_{\gamma_{k+1}}(x_{k+1})$ can be shown to satisfy 
\begin{align*}
V_{\gamma_{k+1}}(x_{k+1}) \leq \varrho \lambda_{\gamma_{k+1}} V_{\gamma_k}(x_k) + \chi(\|w\|_{\mathcal{L}_{\infty}}).
\end{align*}
In a similar manner, assuming that $\|w\|_{\mathcal{L}_{\infty}} \leq M_w$, one has 
\begin{align*}
\mathbb{E}\Big[V(x_{k+1}) \mathbbm{1}_{\{\overline{s}_{k+1}=\overline{s}\}} \Big] &= \sum_{\ell=0}^{1} \mathbb{E} \Big[V_{\gamma_{k+1}}(x_{k+1}) \mathbbm{1}_{\{\overline{s}_{k+1}=\overline{s}\}} \mathbbm{1}_{\{\gamma_{k+1}=\ell\}} \Big] \\
&\leq \varrho \overline{\theta}(\overline{s})\sum_{\overline{s}' \in \overline{S}} \overline{P}(\overline{s}', \overline{s}) \mathbb{E} \Big[V_{\gamma_k}(x_k)\mathbbm{1}_{\{\overline{s}_{k}=\overline{s}' \}} \Big]\\
&\qquad\qquad\qquad+  \chi(M_w) \mathbb{E}\Big[\mathbbm{1}_{\{\overline{s}_{k+1}=\overline{s}\}}\Big].
\end{align*}
Following the same argument as in the proof of Theorem \ref{thm-asas}, the vectorial form of the above inequality can be written as 
\begin{align*}
\overline{V}_{k+1}=\varrho \text{diag}\big(\overline{\theta}(\overline{s}_{1}), \ldots, \overline{\theta}(\overline{s}_{NM})\big) \overline{P} \overline{V}_{k}+\chi(M_w)\overline{\pi}_{k+1}
\end{align*}
with $\overline{V}_{k}=\Big[\mathbb{E}\big[V(x_{k}) \mathbbm{1}_{\{\overline{s}_{k}=\overline{s}\}} \big]\Big]_{\overline{s} \in \overline{S}}$ and $\overline{\pi}_{k+1}=[\mathbb{P}\{\overline{s}_{k+1}=\overline{s}\}]_{\overline{s} \in \overline{S}}$. 
We may then write  
\begin{align*}
\mathbb{E}\big[V(x_{k+1})\big] &\leq \varrho \|\text{diag}\big(\overline{\theta}(\overline{s}_{1}), \ldots, \overline{\theta}(\overline{s}_{NM})\big) \overline{P}\|_{1} \mathbb{E}\big[V(x_{k})\big]\\
&\qquad\qquad\qquad+\chi(M_w).
\end{align*}
Now, if the condition in \eqref{ineq: practical-stable-in-probability} holds with $\eta \in [0, 1)$, then
\begin{align*}
\mathbb{E}\big[V(x_{k+1})\big] &\leq \eta  \mathbb{E}\big[V(x_{k})\big] + \chi(M_w) \\
&\leq \eta^{k+1} \alpha_{2}(|x_0|)+\chi(M_w)\big[1+\eta+\cdots+\eta^{k}\big] \\
& \leq \eta^{k+1} \alpha_{2}(|x_0|)+\chi(M_w)\frac{1-\eta^{k+1}}{1-\eta}
\end{align*}
Since $\alpha_{1}(|x_k|) \leq V_{\gamma_k}(x_k)$ holds for all $k \in \mathbb{Z}_{+}, $ due to property (V1) in Assumption \ref{assumption: multiple-lyapunov}, one knows that $\{x_k \, \vert \, \alpha_{1}(|x_k|) \geq c \}$ is a subset of $\{x_k \,\vert\, V_{\gamma_k} (x_k) \geq c\}$. 
Thus, $\mathbb{P}\{\alpha_{1}(|x_k|) \geq c  \} \leq \mathbb{P}\{V_{\gamma_k} (x_k) \geq c\}, \forall c \geq 0, k \in \mathbb{Z}_+$. 
By Markov's inequality, the probability that the system state $x$ exits a predefined target set $\Omega_{s}=\{x \vert |x| \leq \Delta\}$ can be computed as
\begin{align*}
\lim_{k \rightarrow \infty}\mathbb{P}\big\{|x_k| \geq \Delta + \epsilon \big\} &\leq \lim_{k \rightarrow \infty}\mathbb{P}\big\{V(x_k) \geq \alpha_{1}(\Delta + \epsilon) \big\} \\
&\leq \lim_{k \rightarrow \infty}\frac{\mathbb{E}\big[V(x_k)\big]}{\alpha_{1}(\Delta + \epsilon)} \leq \frac{\chi(M_w)}{(1-\eta)\alpha_{1}(\Delta+\epsilon)}.
\end{align*}
\end{IEEEproof}

\end{document}